\documentclass[11pt]{amsart}  
\usepackage{amsmath,amsfonts,amssymb}
\usepackage{times} 

\newtheorem{lemma}{Lemma}[section]
\newtheorem{theorem}[lemma]{Theorem}
\newtheorem{corollary}[lemma]{Corollary}
\newtheorem{assumption}{Assumption}
\newtheorem{definition}[lemma]{Definition}

\newtheorem{remark}[lemma]{Remark}

\newcommand{\Z}{\mathbb{Z}}
\newcommand{\R}{\mathbb{R}}

\newcommand{\RR}{\mathbb{R}} 
\newcommand{\ip}[2]{{\langle#1,#2\rangle}}
\newcommand{\MD}{{\mathsf D}}
\newcommand{\FD}{{D}}
\newcommand{\MM}{{\mathsf M}}

\newcommand{\N}{{\mathbb N}}

\newcommand{\A}{\mathcal{A}}

\newcommand{\Fc}{\mathcal{F}}

\newcommand{\E}{\mathbb{E}}

\newcommand{\Hbb}{\mathcal{G}}  
\newcommand{\PP}{\mathbb{P}}
\newcommand{\VV}{\mathbb{V}}

\newcommand{\eps}{\varepsilon}

\newcommand{\norm}[1]{\lvert\!\lvert #1 \rvert\!\rvert}
\newcommand{\Norm}[1]{\lvert\!\lvert\!\lvert #1 \rvert\!\rvert\!\rvert}

\newcommand{\bpf}[1][Proof]{\begin{proof}[#1]}
\newcommand{\epf}{\end{proof}}

\newcommand{\osc}{\mathop{\rm osc}}
\newcommand{\spn}{\mathop{\rm span}}
\newcommand{\sch}{\mathop{\rm SCH}}
\newcommand{\Dom}{\mathop{\rm Dom}}
\newcommand{\Extr}{\mathop{\rm Extr}}
\newcommand{\Poly}{\mathop{\mathrm{ Poly}}}

\newcommand{\HH}{\mathbb{H}}

\newcommand{\XX}{\mathbb{X}}
\newcommand{\YY}{\mathbb{Y}}

\newcommand{\D}{\mathbb{D}}

\newcommand{\Ac}{\mathcal{A}}
\newcommand{\Ker}{\mathop{\rm Ker}}
\newcommand{\Lip}{\mathrm{Lip}}
\newcommand{\Hol}{\mathrm{Hol}}

\newcommand{\LipT}{\widetilde \Lip}
\newcommand{\supT}{\widetilde \sup}

\newcommand{\Tt}{t_*}
\newcommand{\TT}{T_*}
\newcommand{\TTT}{T_0}
\newcommand{\Ta}{T_1}
\newcommand{\Tb}{T_2}

\newcommand{\bl}{[}
\newcommand{\br}{]}
\newcommand{\eqdef}{\stackrel{\mbox{\tiny def}}{=}}
\newcommand{\Ha}{\mathrm{H}}
\newcommand{\Hb}{\mathcal{H}}
\newcommand{\Sa}{\mathrm{S}}

\newcommand{\Aa}{\mathrm{A}}
\newcommand{\Ab}{\mathcal{A}}

\renewcommand{\H}{\mathbf{H}}
\newcommand{\V}{\mathbf{V}}
\newcommand{\la}{\langle}
\newcommand{\ra}{\rangle}

\begin{document}
\title{Malliavin Calculus for Infinite-Dimensional
Systems with Additive Noise}
\date{October, 2006}
\author{Yuri Bakhtin}
\address{    School of Mathematics\\
            Georgia Institute of Technology\\
            Atlanta, GA 30332-0160}
\email{bakhtin@math.gatech.edu}
\author{Jonathan C. Mattingly}
\address{   Mathematics Department\\
            Duke University, Box 90320\\
            Durham, NC 27708-0320}
\email{jonm@math.duke.edu}

\begin{abstract}
  We consider an infinite-dimensional dynamical system with polynomial
  nonlinearity and additive noise given by a finite number of Wiener
  processes.  By studying how randomness is spread by the dynamics, we
  develop in this setting a partial counterpart of H\"ormander's
  classical theory of Hypoelliptic operators.  We study the
  distributions of finite-dimensional projections of the solutions and
  give conditions that provide existence and smoothness of densities
  of these distributions with respect to the Lebesgue measure.  We
  also apply our results to concrete SPDEs such as a Stochastic
  Reaction Diffusion Equation and the Stochastic 2D Navier--Stokes
  System.
\end{abstract}

\maketitle
\section{Introduction}\label{sc:Intro}

This paper investigates how randomness is spread by an
infinite-dimensional nonlinear dynamical system forced by a finite
number of independent Wiener processes. The randomness is transferred
by the nonlinearity to degrees of freedom other than those where it is
initially injected. It would be very interesting to obtain precise
information on how the randomness is spread. We will instead show that
some transfer happens almost surely.  Though we are fundamentally
interested in infinite-dimensional systems, we begin with a brief
discussion in finite dimensions.

Consider a stochastic differential equation with additive noise:
\begin{equation}\label{eq:sdeRn}
\left\{  \begin{aligned}
  dx_t &= F_0(x_t) dt + \sum_{k=1}^d F_k dW_k(t)\\
  x_0 &=x\in \RR^m
  \end{aligned}\right.,
\end{equation}
where the $W_k$ are independent standard Brownian Motions, $F_0:\RR^m
\rightarrow \RR^m$ is a bounded analytic function and $F_k \in \RR^m$
is a fixed vector for each $k \in \{1,\cdots,d\}$. Given a function
$u_0:\RR^m \rightarrow \RR$, we can define $u(x,t) =
\mathcal{P}_tu_0(x)\eqdef\E_x u_0(x_t)$. (Here the notation for the
expectation $\E_x$ reinforces the fact that $x_0=x$.) Then $u(x,t)$
solves the backward-Kolmogorov equation $\partial_t u = \mathcal{L} u$
with $u(x,0)=u_0(x)$ and
\begin{equation*}
  \mathcal{L}= F_0\cdot \nabla+ \frac12 \sum_{k=1}^d (F_k \cdot
  \nabla)^2\ . 
\end{equation*}
If the $\spn\{F_1,\cdots,F_d\}= \RR^m$, then the differential operator
is uniformly elliptic. In this case, it is classical that $u(x,t)$ is
a smooth function of $(x,t)$ and that $u(x,t)=\int_{\RR^m}
\rho_t(x,y)u_0(y)dy$, where $\rho$ is a smooth, positive function of
$(t,x,y)$. The function $\rho_t(x,y)$ is called the density of $x_t$
starting from $x_0=x$ (see \cite{bass99DEO}). The fact that $\rho_t$
is smooth and positive is a direct consequence of the randomness
spreading through all of the degrees of freedom.

If  $\text{dim}\spn\{F_1,\cdots,F_d\}< m$, then the preceding
conclusions do not necessarily hold. However, if 
\begin{equation}\label{eq:brackets}
   \spn\Big\{F_i, [F_i,F_{k_1}],  [[F_i,F_{k_1}],F_{k_2}],
   \cdots : i\geq 1, k_j \geq 0, j \geq 1  \Big\}= \RR^m,
\end{equation}
then the system is hypoelliptic and the above conclusions again hold
(save positivity).  Here $[F_j,F_k]=(F_j\cdot\nabla) F_k -
(F_k\cdot\nabla) F_j$ is the Lie bracket (or commutator) of the two
vector fields. Since in our setting the $F_j$ are constant for $j\geq
1$, only brackets with $F_0$ produce non-zero results.  The fact that
the system is hypoelliptic follows from H\"ormander's pioneering work.
In particular, it falls under a generalization of his ``sum of
squares'' theorem.  (The principal part of $\mathcal{L}$ is the sum of
squares of vector fields.)

In the 1970's and 1980's there was a large body of work to develop a
probabilistic understanding of H\"ormander's theorem and related
concepts by looking directly at \eqref{eq:sdeRn} rather than the PDE
for $u_t(x,y)$ and $\rho_t(x,y)$. This line of work was initiated by
Malliavin and contained substantial contributions from Bismut,
Stroock, Kusuoka, Norris and others. The tools developed to address
this question go under the heading of Malliavin Calculus (see
\cite{b:Malliavin78,b:KusuokaStroock84,b:Bell87}) which might well be
described as the stochastic calculus of variations.

We are interested in developing a version of these results in infinite
dimensions ($m=\infty$). We wish to understand which of the previous
conclusions hold if we assume that some variation on
\eqref{eq:brackets} holds with $m=\infty$, where the SDE in
\eqref{eq:sdeRn} is replaced by a stochastic partial differential
equation (SPDE). From the beginning, it is clear that we cannot work
directly with the density $\rho_t(x,y)$ since in infinite dimensions
there is no Lebesgue measure. Ideally, we would like to find a natural
replacement for Lebesgue measure in the setting of a given equation.
For the moment this escapes us, so we will make statements about the
finite-dimensional projection of $\rho$ and the spectrum of the
Malliavin covariance matrix (see Section \ref{secMalliavin}). One
might reasonably ask if we could ever expect some form of
H\"ormander's condition to hold when the dimension $m$ is infinite but
the number of Brownian forcing terms $d$ is finite. In
\cite{short:EMattingly01,b:Romito04}, it was shown that the 
finite-dimensional Galerkin truncations of the Navier-Stokes equations
satisfy H\"ormander's condition for hypoellipticity independent of the
order of the truncation. And thus in some sense,  H\"ormander's
condition holds, at least formally, for the whole SPDE ($m=\infty$).

In \cite{BaudoinTeichmann05HID}, the authors treat the case when the
infinite-dimensional ($m=\infty$) evolution generates a fully
invertible flow and prove conditions guaranteeing the existence of a
density of the finite-dimensional marginals. Because they assume the
dynamics generate a flow, their exposition more closely mirrors the
finite-dimensional treatment. In particular, they are able to handle
diffusion constants which depend on the state of the process. We will
see that our treatment will lead to objects not adapted to the Wiener
filtration, which makes the general diffusion case more difficult.

Because our PDEs generate only a semi-flow (and not a full flow), we
cannot apply directly the same proofs developed using Malliavin
Calculus in the finite-dimensional setting or the infinite-dimensional
extensions given in \cite{BaudoinTeichmann05HID}.  However we will see
that we can modify the proofs to produce the desired results. D. Ocone
\cite{b:Ocone88} first used related ideas in the infinite-dimensional
case when the equations were linear in the solution and the noise; and
hence, an explicit formula exists for the solution. In
\cite{MattinglyPardoux05}, the 2D Navier-Stokes equations are
considered with additive noise. The techniques used there are very
close to those used here. However, there the scope is more limited.
The calculations are done in coordinates which leads to the
restriction that the forcing is diagonal in the same basis.  In both
cases, as in \cite{b:Ocone88}, the time reversed adjoint of the linear
flow is used to propagate information backwards in time. This leads to
a need for estimates on Wiener polynomials with non-adapted
coefficients.  In \cite{MattinglyPardoux05}, only second-order
polynomials were considered. Here, by simplifying and streamlining the
proofs, we can handle general polynomials of finite order. This allows
us to treat PDEs with more general polynomial nonlinearities. Lastly,
we observe that if one is only interested in the existence of a
density, one can jettison over two-thirds of the paper and all of the
technically involved sections.

To further motivate this article, we mention that the type of
quantitative estimates on the spectra of the Malliavin covariance
matrix obtained in this paper is a critical ingredient in the recent
proof of unique ergodicity of the two-dimensional Navier Stokes equations
under the type of finite-dimensional forcing considered in this note (
see \cite{b:HairerMattingly04}). The results of this paper are a major
step towards proving similar results for other SPDEs.

In \cite{EckmannHairer00UIM}, the ergodicity of a degenerately forced
SPDE was also proven using techniques from Malliavin calculus. In
contrast to our setting, there infinitely many directions were forced
stochastically.  However, the  structure of the forcing was
such that it caused the asymptotic behavior for the high spatial modes
to be close to that of an associated linear equation. The type of
analysis used there does not seem to be possible in our setting.

Independently, and contemporaneously to this work M. Wu completed a thesis
\cite{Wu2006Thesis} which carried out the program from
\cite{MattinglyPardoux05,b:HairerMattingly04} to prove the unique
ergodicity of a degenerately forced Boussinesq equation. Since this
equation has a quadratic nonlinearity, he was able to use the
technical lemmas from \cite{MattinglyPardoux05}. However, he also
proved a more general technical lemma which can be used to prove the
existence, but not smoothness, of finite-dimensional marginal
distributions. He used this result to prove the existence of 
finite-dimensional marginal densities for a degenerately forced cubic
reaction-diffusion equation. The technical lemma is similar to
Lemma 5.1 and Proposition 5.2 from \cite{MattinglyPardoux05},
though the proof given is slightly different. In this note, we have
used the other, though related, approach from
\cite{MattinglyPardoux05}, since (at least for us) it is more
straightforward to use it to obtain the quantitative estimates needed to
prove smoothness.

While this paper was in its final stages of completion, the authors
became aware of a recent preprint
\cite{AgrachevKuksinSarychevShirikyan06FDP} where it is proven that
finite-dimensional projections of a randomly forced PDE's Markov
transition kernel are absolutely continuous with respect to Lebesgue
measure if a certain controllability condition is satisfied.  While
connections between controllability and the existence of densities are
not surprising given what is known for maps and SDEs (see
\cite{b:Kliemann87,b:BenArousLeandre91} for example), the strength of
the results in \cite{AgrachevKuksinSarychevShirikyan06FDP} is that
they do not require the forcing to be Gaussian.  They only need that
it satisfies a more general condition of decomposability. However,
that approach presently does not provide smoothness of the densities.

\vspace{1ex}
\noindent \textbf{Organization:} In Section \ref{sc:framework}, we
introduce the abstract setting for the rest of the paper. In
Section~\ref{sec:basicResults}, we give the main results of the paper
in a simplified form which is sufficient for the applications we
present. Principally, we give results ensuring the existence and
smoothness of the finite-dimensional projections of the Markov
transition kernels. In Section \ref{sec:Applications}, we specialize
the abstract framework and apply it to a scalar reaction-diffusion
equation and the two-dimensional Navier--Stokes equation. In
Section~\ref{secMalliavin}, we give a brief introduction to the ideas
from Malliavin calculus we need.  In Sections \ref{sec:GeneralResults}
and \ref{sc:genResults}, we respectively state and prove the principal
results in their full generality. In Section \ref{sec:refinments}, we
give a number of generalizations and refinements tailored to the needs
of the arguments in \cite{b:HairerMattingly04} which prove the unique
ergodicity of the system as already mentioned. The estimates on the\
spectrum of the Malliavin Covariance matrix in this paper constitute
one of the principal ingredients of that work.  In
Section~\ref{sc:non-adapted_stuff}, we give the necessary abstract
results on non-adapted polynomials of Wiener processes, one of the
main technical tools of the paper. In the remaining two sections, we
give a number of auxiliary lemmas needed in the proofs.

\vspace{1ex}
\noindent\textbf{ Acknowledgments: } This work grew from a joint work
of JCM with \'Etienne Pardoux whom he thanks for the many fruitful,
interesting and educational discussions which laid the ground work for
this work.  YB thanks the hospitality of Duke University during the
academic year 2004--2005 when the bulk of this work was done.  The
authors also thank Trevis Litherland, Scott McKinley, and Natesh Pillai for reading and
commenting on a preliminary version of this paper. JCM was supported
in part by the Sloan Foundation and by an NSF PECASE award
DMS04-49910.


\section{General Setting}\label{sc:framework}

In this section we introduce the framework to define and study
solutions of a stochastic evolution equation in a Hilbert space
\begin{equation}
\label{eq:spde}\left\{
  \begin{aligned}
    du(t)=&L(u) dt + N(u)dt+f(t)dt+\sum_{k=1}^dg_kdW_k(t),\\
    u(0)=&u_0.
  \end{aligned}\right.
\end{equation}

The three components of the framework are: the space where solutions
are to be defined; the deterministic part of the r.h.s. (the
drift), namely, the autonomous part given by the vector field $L(u) +
N(u)$ and the non-autonomous part $f$; and the noise $dW$. The operator $L$ is
assumed to be linear while $N$ contains all of the nonlinear
terms. More details are given in the following.

The first component of the framework is the space where the solutions
are going to live.  We need two separable Hilbert spaces
$\HH$ and $\VV$, with norms $|\cdot|$ and $\|\cdot\|$ generated by
inner products $\langle\cdot,\cdot\rangle$ and
$\langle\!\langle\cdot,\cdot\rangle\!\rangle$, respectively.  We assume
that $\VV$ is compactly embedded and dense in $\HH$ so that $\HH$ is
compactly embedded and dense in $\VV'$, the dual of $\VV$. Hence
$\ip{\cdot}{\cdot}$ also gives the duality pairing between $\VV$ and
$\VV'$. We also assume that $|v|\le\|v\|$ for any $v\in\VV$.

We shall assume that there is a set $\HH_0 \subset \HH$ such that with
probability one
\begin{equation}
\label{eq:where_u_lives}
u\in C\big([0,T],\HH\big)\cap  C\big((0,T],\VV\big)
\end{equation}
for  all $u_0\in \HH_0$.

The deterministic external force $f$ is a bounded $\HH$-valued function
defined on a time interval $[0,T]$.
$L$ is a linear operator with values in $\VV'$ defined on a subspace
of $\HH$. The restriction of $L$ to $\VV$ is a continuous operator
$\VV\to\VV'$.

The nonlinear vector field $N:\VV\to \HH$ will be assumed to be a
continuous polynomial (defined below) with zero linear and constant
part.  It is convenient to introduce the notation
\begin{equation*}
  u^{\otimes j}=(\underbrace{u,\ldots,u}_{j}).
\end{equation*}
Often, for a function $Q$ of $j$ variables, we shall write  
\begin{equation*}
  Q(u)=Q(u^{\otimes j})=Q(\underbrace{u,\ldots,u}_{j}).
\end{equation*}

\begin{definition}
  Given two Banach spaces $\XX$ and $\YY$, we say that $F:\XX\to \YY$
  is a {\it continuous polynomial} of positive integer degree~$m$ if
\begin{equation*}
  F(x)=F(x^{\otimes m})
\end{equation*}
for some map $F:\XX^m\to\YY$ such that
\begin{equation*}
  F(x_1,\ldots,x_m)=F_0+F_1(x_1)+F_2(x_1,x_2) + \ldots +
  F_m(x_1,\ldots,x_m),
\end{equation*}
where all functions $F_j:\XX^j\to\YY$ are multilinear (i.e., linear in
each variable), symmetric (i.e., invariant under argument permutations),
and continuous.
\end{definition}

 Hence our assumption on $N$ states that
\begin{align*}
  N(u_1,\ldots,u_m) = N_2(u_1,u_2) + \ldots +
N_m(u_1,\ldots,u_m),
\end{align*}
where all functions $N_j:\VV^j\to\HH$ are multilinear, continuous, and
symmetric. For notational convenience we will write $F(u)=L(u)+N(u)$.

Finally, our probability space is $(\Omega,\Fc,\PP)$, where
$\Omega=C([0,T],\R^d)$, and $\PP$ is the standard Wiener measure on
$\Omega$ equipped with the completion~$\Fc$ of the Borel
$\sigma$-algebra induced by the $\sup$-norm.  The noise $W$ is given
by the canonical map $W(\omega)=\omega$, $\omega\in\Omega$. The $g_i$
from equation \eqref{eq:spde} are fixed elements of $\Dom(L)\cap \VV$.
(Here and in the sequel, we use the notation $\Dom(L)=\{v\in\HH:
L(v)\in\HH\}$.) Letting $e_1,\ldots,e_d$ denote the standard basis in
$\R^d$, we define a linear map $G:\R^d\rightarrow \Dom(L)\cap \VV$
by $g_1=G e_1,\ldots,g_d=G e_d $.

As usual, the stochastic equation~\eqref{eq:spde} is simply shorthand
for the following integral equation:
\begin{equation}
\label{eq:spde-integral}
u(t,\omega)=u_0+\int_0^tF(u(s,\omega))\,ds+\int_0^tf(s)\,ds+GW(t,\omega).
\end{equation}

We shall always assume that there exists a stochastic semiflow
associated with this equation.  More precisely, we assume that there
is a family of operators
\begin{equation*}
  \Phi_t:C([0,t],\R^d)\to \VV,\quad t\in(0,T],  
\end{equation*}
such that if $u(t)=\Phi_t(W[0,t])$ for all $t\in[0,T]$ with
probability~1,  then $u$ is a solution of \eqref{eq:spde-integral}
satisfying~\eqref{eq:where_u_lives}. Here $W[0,t]$ is the restriction
of~$W$ to $[0,t]$. We stress that though the initial data $u_0$ may be
in $\HH\setminus\VV$, the solution is assumed to be in
$C((0,T],\VV)$.

\section{Basic Results}
\label{sec:basicResults}
Our main results are the absolute continuity of the distribution of the
projection of $\Phi_T(W)$ on a finite-dimensional space with respect
to the Lebesgue measure on that space, and the smoothness of the density.

We will need some conditions on the linearization of the
system~\eqref{eq:spde}.  Let $J_{s,t}:\HH\to\VV,0<s\leq t$ solve the
equation in variations:
\begin{equation}\label{eq:J}
  \left\{  \begin{aligned}
      \frac{\partial}{\partial t}J_{s,t}\phi=&(\FD
      F)(u(t))J_{s,t}\phi, \quad s< t,\\ 
      J_{s,s}\phi=&\phi,\quad \phi\in\HH.  
  \end{aligned}\right.
\end{equation}
Here $(\FD F)(x)h$ is the Fr\'echet derivative of $F$ at a point
$x\in\VV$ applied to a tangent space vector $h\in\VV$. Hence, for each
$x$ we have $(\FD F)(x):\VV\to\VV'$.  Notice that Fr\'echet derivatives
of $F$ of all orders are well-defined (see
Lemma~\ref{lm:higher_derivatives_of_multilinear}.)

\begin{assumption}\label{a:J1}
  With probability one, there is a unique solution $J_{s,t}\phi$ to
  the equation~\eqref{eq:J} for every
  $\phi\in\HH$  and $s \in (0,T)$ with
\begin{equation}
  \label{eq:where_J_lives}
  J_{s,t}\phi\in C\big([s,T],\HH\big)\cap L^2\big([s,T],\VV\big),\quad \frac{\partial}{\partial
    t}J_{s,t}\phi\in L^2\big([s,T],\VV'\big), 
\end{equation}
where $J_{s,t}$ and $\frac{\partial}{\partial t}J_{s,t}$ are
considered as functions of $t$.  
\end{assumption}

We are also going to consider the time reversed adjoint of $J_{s,t}$
denoted by $K_{s,t}:\HH\to\VV,s\leq t$ and defined by the backward
equations
\begin{equation}\label{eq:adjoint_equation}
  \left\{  \begin{aligned} 
      \frac{\partial}{\partial s}K_{s,t}
      \phi=&-(DF)^*(u(s))K_{s,t}\phi,\quad s<t,\\ 
      K_{t,t}\phi=&\phi,\quad  \phi\in\HH.
    \end{aligned}\right.
\end{equation}
Here $(DF)^*(y):\VV\to\VV'$ is the adjoint operator for $(DF)(y)$.
(We identify $\VV$ and $\VV''$.) 
\begin{assumption}\label{a:K1}
  With probability one, for any $0<t_0<t\leq T$ and $\phi \in \HH$,
  there is a unique solution $K_{s,t}\phi$ of
  equation~\eqref{eq:adjoint_equation} that satisfies
\begin{equation}
  \label{eq:where_K_lives}
  K_{s,t}\phi\in C\big([t_0,t],\HH\big) \cap L^2\big([t_0,t],\VV\big),\quad
  \frac{\partial}{\partial s}K_{s,t}\phi\in L^2\big([t_0,t],\VV'\big), 
\end{equation}
where $K_{s,t}$ and $\frac{\partial}{\partial s}K_{s,t}$ are
considered as functions of $s$. 
\end{assumption}

 \subsection{Existence of densities}
\label{sec:ExtDenBasic}
 To begin understanding how the randomness spreads through the phase
 space, we now introduce an increasing collection of sets which
 characterize some of the directions excited. In Section
 \ref{sec:GeneralResults} we will give a more completed, though more
 complicated, description of the directions excited. However, for many
 cases, the results of this section are sufficient.

For any positive integer $n$, we introduce the subset $\Hbb_n$ of $\VV$
defined recursively as follows.  For $n=1$, we set
$\Hbb_1=\spn\{g_1,\ldots,g_d\}$.   For $n>1$, $\Hbb_n$ is defined via
 $\Hbb_{n-1}$: we set
\begin{multline}\label{def:G}
  \Hbb_n\eqdef\spn\Bigl(
  \Hbb_{n-1}\bigcup\Big\{N_m(g,g_{k_1},\ldots,g_{k_{m-1}}) \in\VV : \\g
    \in\Hbb_{n-1}\cap\VV\cap \Dom(L), k_i \in\ \{1,\cdots,d\}\Big\}\Bigr). 
\end{multline}

Finally, we introduce $\Hbb_\infty=\spn\big(\bigcup_n\Hbb_n\big)$. The following
result is a specialization of Theorem \ref{th:absolute_continuityB},
which is given later.

\begin{assumption}\label{a:D2}
  For all $u_0 \in \HH_0$, there is a constant $J^*(T,u_0)$ such that
  \begin{equation}
    \label{eq:assumption_on_2nd_moment_of_J}
    \sup_k \sup_{0 < s<t\leq T}\E|J_{s,t}g_k|^2\leq J^*(T,u_0).
  \end{equation}
\end{assumption}

\begin{theorem}\label{th:absolute_continuityA} Assume that Assumptions
  \ref{a:J1}, \ref{a:K1}, and  \ref{a:D2} hold.  Suppose that $S \subset
  \HH$ is a finite-dimensional linear subspace in $\Hbb_\infty$.  Then
  the distribution of the orthogonal projection $\Pi_S u(T)$ on $S$ is
  absolutely continuous with respect to the Lebesgue measure on~$S$.
\end{theorem}

\subsection{Smoothness of densities}
\label{sec:smoothDenBasic} To prove smoothness of the density obtained
in Theorem~\ref{th:absolute_continuityA}, we need stronger assumptions
than those made in that theorem.  We will replace the assumption of
continuity and finiteness of the first derivative in time of $K$ and $u$ with
assumptions on the moments of the Lipschitz coefficients in time of related
processes.

\begin{assumption}\label{a:solReg} In addition to the standing assumptions from
  section \ref{sc:framework}, the following conditions hold: 
 For every $u_0 \in \HH_0$ there exists a fixed $\TTT \in[0,T)$ and
    constants $u^*_p(\TTT,T,u_0)$, for all integers $p \geq 1$, so that
    \begin{align}\label{eq:assumption_moment}
      \E \sup_{\TTT \leq t\leq T} \|X(t)\|^p &\leq  u^*_p(\TTT,T,u_0),\\
      \label{eq:assumption_momentLip}
       \E \sup_{\TTT\leq s<t\leq T}
       \left[\frac{|X(t)-X(s)|}{|t-s|}\right]^p &\leq 
      u^*_p(\TTT,T,u_0), 
    \end{align}
where $X(t)=u(t)-GW(t)$.
\end{assumption}

\begin{assumption}\label{a:K2a} In addition to Assumptions  \ref{a:J1}
  and \ref{a:K1} there exists a $\TTT$ so that for every $u_0 \in \HH_0$
  and $p \geq 1$ there exists a constant $K^*_p(\TTT,T,u_0)$ with
    \begin{align}
      \label{eq:assumption_on_moments_of_K}
      \E\left[\sup_{\TTT\leq s<t\leq
          T}|K_{s,t}|^p_{\VV\to\VV}\right]\leq
      K^*_p(\TTT,T,u_0),\\
      \E\left[\sup_{\TTT\leq s<t\leq
          T}\left(\frac{|K_{s,T}-K_{t,T}|_{\VV\to\HH}}{|t-s|}\right)^p\right]\leq
      K^*_p(\TTT,T,u_0),    \label{eq:assumption_on_Lip_of_K}
    \end{align}
    where $|\cdot|_{\VV\to\VV}$ denotes the norm of a linear operator mapping $\VV$
    to itself and $|\cdot|_{\VV\to\HH}$ from $\VV$ into $\HH$.
\end{assumption}

\begin{assumption}\label{a:J2a}There exists an
  $\alpha \in [0,1)$ such that for any $u_0 \in \HH_0$ and $p \geq1$
  there is a constant $D^*_p(T,u_0)$ such that
    \begin{align}\label{eq:JmomentWSimple}
      \sup_k \;\E\sup_{0\leq s<t \leq
        T}\|J_{s,t}g_k\|^p &\leq D^*_p(T,u_0),\\
      \E\sup_{0\leq s<t \leq T}\Bigl[(t-s)^\alpha|J_{s,t}|_{\HH \to
          \VV}\Bigr]^p &\leq D^*_p(T,u_0)\;.
      \label{eq:pth_momment_JSmooth}
    \end{align}
 \end{assumption}

 Lastly, we need the following definition which further restricts the
 class of polynomial nonlinearities we will treat.
 \begin{definition}\label{poly1}
   We define $\Poly_1(\VV,\HH)$ as the set of continuous polynomials $Q:\VV
   \to \HH$, with $Q=\sum_{i=1}^k Q_i$ for some $k$, where the $Q_i$
   are homogeneous $i$-linear terms satisfying the following bound for
   some $C_i$ and all $u_i \in \VV$:
   \begin{equation}
     \label{eq:Poly1}
     |Q_i(u_1,\cdots,u_i)|_{\VV'} \leq C_i |u_1|\|u_2\|\cdots\|u_i\|.
   \end{equation}
 \end{definition}

 We are now in a position to state precisely our first result on the
 smoothness of the projections of transition densities.  
\begin{theorem}\label{th:smoothnessA}
  In the setting of Section \ref{sc:framework}, assume that
  Assumptions \ref{a:solReg}, \ref{a:K2a} and \ref{a:J2a} hold.  Let
  $S \subset \VV$ be a finite-dimensional subspace of $\Hbb_n$ for some $n$.
  If $N$ is a continuous polynomial in $\Poly_1(\VV,\HH)$ then the density of
  $\Pi_S u(T)$ with respect to Lebesgue measure exists and is a
  $C^\infty$-function on $S$.
\end{theorem}

\section{Applications}\label{sec:Applications}

We now specialize our setting, restricting ourselves to the case where the linear
operator $L$ is dissipative and dominates the nonlinearity. At the end
of the section we will fit a reaction-diffusion equation and the 2D
Navier-Stokes equation into the setting we now describe.

Let $L$ be a positive, self-adjoint linear operator on a Hilbert space
$\H$. Additionally, assume that $L$ has compact resolvent.  Hence $L$
has a complete orthonormal eigenbasis $\{e_k : k=1,2,\cdots\}$, with
real eigenvalues $0<\lambda_1 \leq \lambda_2 \leq \cdots$ such that
$\lim \lambda_k = \infty$. For $s \in \R$, we define the inner product
\begin{equation*}
   \la u,v \ra_s =\sum_{k=1}^\infty \lambda_k^{2s} \la u,e_k\ra  \la e_k,v\ra 
\end{equation*}
and the norm $|u|^2_s = \la u, u \ra_s$. We define the spaces
\begin{equation*}
  \V^s= \{ u  \in \H : |u|_s < \infty\}
\end{equation*}
and observe that $\V^{-1}$ is the dual in $\H$ of $\V^1$ and that $L$
maps $\V^1$ to $\V^{-1}$ and $\V^2$ to $\H$.  We assume that $N \in
\Poly_1(\V^1,\H)$ and that $f:[0,T]\rightarrow \V^1$ is uniformly
bounded. (See definition \ref{poly1}.)
\begin{lemma}\label{l:ExistSolExmpl}
  In the setting above, for any $u_0 \in \H$ equation \eqref{eq:spde}
  has a unique strong solution $u$, generated by a stochastic flow
  $\Phi_t: \H \rightarrow \V^2$, satisfying
  \begin{equation*}
    u \in C\big([0,T] ; \H\big) \cap C\big((0,T]; \V^2\big).
  \end{equation*}
  In addition, if $\phi \in \H$, then there exists a unique solution $J_{s,t}
  \phi$ to equation \eqref{eq:J}, for all $0 < s\leq t\leq
  T$. Furthermore
  \begin{equation*}
    J_{s,t}\phi \in C\big([s,T]; \H\big) \cap  C\big((s,T]; \V^2\big)
    \cap  L^2\big([s,T]; \V^1\big)  
  \end{equation*}
as a function of $t$.

Lastly, if $\phi \in \H$, then there exists a unique solution $K_{s,t}
  \phi$ to equation \eqref{eq:adjoint_equation}, for all $0 <t_0 < s\leq t\leq
  T$. Furthermore,
  \begin{equation*}
    K_{s,t}\phi \in C([t_0,t]; \H) \cap  C([t_0,t]; \V^2) \cap  L^2([t_0,t]; \V^1) 
  \end{equation*}
as a function of $s$.
\end{lemma}

\begin{proof}
  Most of the results follow from results about deterministic, time
  inhomogeneous equations found in \cite{SellYou02DEE}. As is often
  done (see for example \cite{b:Fl94,b:DaZa96}), we begin by setting
  $u(t)=X(t) + GW(t)$. Then $X(t)$ satisfies a standard PDE
  \begin{equation*}
     \frac{\partial}{\partial t}X(t) = F(X(t),t),
  \end{equation*}
  where the random right hand side is given by $F(x,t)= L(x) +
  N(x+GW(t))$. Once it is demonstrated that this equation has a unique
  solution for every $u_0 \in \H$ and almost every $W$, we have
  constructed a stochastic flow, since all initial conditions can share a
  single exceptional set in the probability space.  Clearly,
  $(x,t)\mapsto N(x+GW(t))$ is
  almost surely uniformly bounded in $\H$ on $\{ (x,t) : |x|_1 + t <
  C\}$, for all $C>0$. Furthermore, it is a H\"older continuous function
  of time for all H\"older exponents less than $1/2$. The quoted
  existence, uniqueness, and regularity for $u$ then follows from Lemma
  47.2 from \cite{SellYou02DEE} applied to the above equation for $X(t)$.

  All of the quoted results on the linearization except the fact that the
  solution is in $L^2([s,T]; \V^1)$, follow from Theorem 49.1 from
  \cite{SellYou02DEE} by arguments similar to those just employed.
    To see that the solutions are in $L^2([s,T]; \V^1)$, take the inner
  product with $v(t)=J_{s,t} \phi$  to obtain
  \begin{align}
     \frac{\partial}{\partial t} |v(t)|^2_0 &=  - |v(t)|^2_1 +
     \ip{DN(u(t))v(t)}{v(t)}_0\notag \\
     &\leq   -\frac12 |v(t)|^2_1 +C(1+  |u(t)|_1^{2(m-1)}|v(t)|^2_0) \;.
  \end{align}
Since this implies that
\begin{align*}
 \int_s^T |v(t)|_1^2 dt \leq C\Bigl(|\phi|^2_0 + \sup_{t \in[s,T]}  |u(t)|_1^{2(m-1)} |v(t)|^2_0\Bigr) < \infty,
\end{align*}
we are done.

The proofs of the statements for the adjoint linearization are the same
as for the linearization after one observes that since  $N \in
\Poly_1(\V^1,\H)$ we have
\begin{align*}
  \sup_{|v_i|_0=1} \ip{DN_j^*(u)v_1}{v_2} &= \sup_{|v_i|_0=1}
  \ip{DN_j(u)v_2}{v_1} \leq  |N_j(v_2,u,\cdots,u)|_0|v_1|_0 \\&\leq
  C|v_2|_0|v_1|_0  |u|_1^{j-1} \leq  C|u|_1^{j-1}.
\end{align*}
\end{proof}

\begin{corollary}\label{c:A1and2OK}
  Setting $\VV=\V^1$, $\HH=\H$, $|\cdot|=\lambda_1|\cdot|_0$, and
  $\|\cdot\|=|\cdot|_1$, the standing assumptions of
  Section \ref{sc:framework} and Assumptions \ref{a:J1} and \ref{a:K1}
  hold in the setting of this section.
\end{corollary}

\subsection{A Reaction-Diffusion Equation}

Consider the following reaction-diffusion equation
\begin{equation}
  \label{eq:RDeqn}
  \left\{\begin{aligned}
    du(x,t) &= \big[\nu \Delta u(x,t) +
    N\big(u(x,t)\big) \big]dt + \sum_{k=1}^dg_k(x) dW_k(t), \\
    u(x,0)&=u_0(x),
  \end{aligned}\right.
  \end{equation}
with $x \in [0,1]$,
\begin{equation*}
  N(u)= \sum_{k=0}^{2q+1} a_k u^k,
\end{equation*}
with $a_k \in \RR$ and $a_{2q+1} < 0$, and under the Dirichlet
boundary conditions
\begin{equation}
  \label{eq:BC}
  u(t,0)=u(t,1)=0 \qquad\text{for all $t \geq 0$}.
\end{equation}
Since in one dimension there exists a constant $C$ so that for any $f
\in \V^1$, $|f|_\infty \leq C |f|_1$ where $|f|_\infty$ is the sup-norm, we see that $N$ is a continuous
polynomial  from $\V^1$ to $\H$ and that 
\begin{equation*}
  DN(u)v = DN^*(u)v =\sum_{k=1}^{2q+1} k a_k u^{k-1}v.
\end{equation*}
The following calculation shows that $N \in  {\Poly}_1(\V^1,\H)$. Let
$u_i \in \V^1$ and observe that 
\begin{align*}
  |N_j(u_1,\cdots,u_j)|_0=&\sup_{\substack{v \in \H\\|v|_0=1}} \ip{N_j(u_1,\cdots,u_j)}{v}
  = a_j\sup_{\substack{v \in \H\\|v|_0=1}} \int u_1(x) \cdots
  u_j(x) v(x) dx  \\&\leq|a_j| |u_2|_\infty \cdots  |u_j|_\infty
  |u_1|_0 \leq C  |u_1|_0|u_2|_1 \cdots  |u_j|_1,
\end{align*}
where the last inequality follows from the Sobolev inequality $|u|_\infty
\leq C |u|_1$.

At the end of the example we will address  necessary conditions for the
system to be formally H\"ormander.  For now we address the
technical conditions needed to apply Theorems
\ref{th:absolute_continuityA} and \ref{th:smoothnessA}. In light of
Corollary \ref{c:A1and2OK}, to apply Theorem
\ref{th:absolute_continuityA} we need to verify
 Assumption \ref{a:D2}. Letting $v(t)=J_{s,t}v_0$, we have 
 \begin{align*}
   \frac12 \frac{\partial }{\partial t} |v(t)|^2_0 &= - \nu |v(t)|^2_1 +
   \ip{DN(u(t))v(t)}{v(t)}_0\\
   &\leq -\nu |v(t)|^2_0 + K_1  |v(t)|^2_0,
 \end{align*}
since $\sup_{a \in \RR} N'(a) \leq K_1/\lambda_1$ for some
$K_1$. Gronwall's inequality then implies
\begin{equation*}
  \sup_{t \in [0,T]}|v(t)|^2_0 \leq  |v_0|^2_0 e^{2(K_1-\nu)T},
\end{equation*}
which translates to 
\begin{equation}\label{eq:JHnormRD}
  \sup_{0\leq s < t \leq T} |J_{s,t}|_{\H \rightarrow \H}^p  \leq
  e^{2p (K_1-\nu)T},
\end{equation}
for all $p >0$. This ensures that Assumption \ref{a:D2} holds.
Having verified all of the assumptions of Theorem
\ref{th:absolute_continuityA}, we have the following result:
\begin{theorem}
  If $g_k \in \V^2$, then  the conclusions of Theorem
\ref{th:absolute_continuityA} hold for equation~\eqref{eq:RDeqn}.
\end{theorem}

We now investigate conditions which guarantee that $\Hbb_\infty$ is dense in $\HH=\H$.
Let $\mathcal{I}_0$ be a finite collection of functions in $\V^2$. Let
$\mathcal{I}$ be the multiplicative algebra generated by
$\mathcal{I}_0$. 

\begin{lemma}
  If $\mathcal{I}$ is dense in $\V^2$, then to ensure that $\Hbb_\infty=\H$ it is sufficient that
  \begin{equation}
    \label{eq:denseCondition}
    \{ f_{1}\cdots f_{k} :  f_j \in \mathcal{I}_0, 1 \leq k \leq 2q \}\subset\Hbb_1
    =\spn\{ g_j : j=1,\cdots,d\}.
   \end{equation}
\end{lemma}

\begin{remark}
  For $\mathcal{I}$ to be dense in $\V^2$, it is sufficient,
  by Stone--Weierstrass, that $\mathcal{I}$ separates points in
  $\V^2$ and if $f(x)=0$ for all $f \in \mathcal{I}$ then
  $x=0$ or $1$.
\end{remark}

We now turn to proving that the density is smooth. In the sequel, we
are going to restrict ourselves to initial data in $\H_0 = \H \cap
C([0,1])$.  It is well known (See \cite{Cerrai99SPT} Proposition 3.2
or \cite{EckmannHairer00UIM}) that for all $p \geq 1$ and $u_0 \in \H_0$
\begin{equation}
  \label{eq:engBoundRD}
  \E  |u(t)|^{p}_0 \leq C(p,t) ( 1+ |u_0|^{p}_0) 
\end{equation}
and
\begin{equation}\label{eq:infBoundRD}
  \E \sup_{t \in[0,T]} \sup_{x \in [0,1]} |u(t,x)|^p \leq C(T,p,u_0)< \infty\ .   
\end{equation}

\subsubsection{Verification of Assumption \ref{a:solReg}} 
Since for any  $T_0 \in(0,T]$ and $t \in [T_0,T]$,  
\begin{equation*}
  X(t)=e^{\Delta(t-T_0)} \Big[ u(T_0)-GW(T_0) \Big] + \int_{T_0}^t
  e^{\Delta(t-s)} \Big[N(u(s))+ \Delta G W(s)\Big] ds  
  \end{equation*}
we know that 
\begin{align*}
  |X(t)|_1 \leq& C|u_{0}|_0 + C\int_{0}^t \frac{e^{-a(t-s)}}{\sqrt{t-s}}
 \Big[ |N(u(s))|_0  + |W(s)|_{\RR^d}\Big]ds\\
  \leq& C\Big( 1+ |u_{0}|_0 +\sup_{s \in [0,T]} |u(s)|^{2q+1}_\infty+ \sup_{s \in [0,T]}|W(s)|_{\RR^d} \Big)
\end{align*}
for some positive $a$ and $C$, and all $t \in [0,T]$.
Applying \eqref{eq:engBoundRD}, \eqref{eq:infBoundRD} and standard
bounds on $\sup |W(t)|$, proves the first part of Assumption
\ref{a:solReg} for any $u_0 \in \H_0$.

Similarly for $0 \leq s < t \leq
T$,
\begin{multline*}
  X(t) - X(s) = e^{\Delta s}\left( e^{\Delta (t-s)} - I\right)
  u_0 \\+ \int_0^s
  e^{\Delta(s-r)} \left[e^{\Delta(t-s)}-I \right]\Bigl[ N(u(r)) +
    \Delta G W(r))\Bigr] dr,
\end{multline*}
which implies there exists a constant $C$ depending on $T$ such that
\begin{multline*}
  | X(t) - X(s)|_0 \leq C|t-s|\Big(1+ |u_0|_0+\ \sup_{r \in[0,T]}|u(r)|_\infty^{2q+1}+\sup_{r \in [0,T]}|W(r)|_{\RR^d}\Big).
\end{multline*}
Again combining standard estimates with \eqref{eq:engBoundRD} and
\eqref{eq:infBoundRD}, we see that the estimate in
\eqref{eq:assumption_momentLip} holds.

\subsubsection{Verification of Assumption \ref{a:K2a}}
Again  setting $v(t)=J_{s,t}v_0$ for $s \in [0,T)$, we have
\begin{align*}
  \frac12 \frac{\partial }{\partial t}|v(t)|_1^2 &= - \nu |v(t)|_2^2 +
  \ip{N'(u(t)) v(t)}{ \Delta v(t)}_0\\
  & \leq  - \nu |v(t)|_2^2 + C\big(1+
  |u(t)|^{2q}_\infty\big)|v(t)|_0|v(t)|_2\\
  & \leq   \frac{C'}{\nu}\big(1+  |u(t)|^{4q}_\infty\big)|v(t)|^2_0,
\end{align*}
which when combined with \eqref{eq:JHnormRD} implies that
\begin{align*}
   \sup_{s\leq t \leq T}|v(t)|_1^2 \leq&  |v(s)|_1^2 + \frac{C'}{\nu} \int_{0}^T
(1+  |u(r)|^{4q}_\infty) |v(r)|^2_0 dr\\
\leq&|v(s)|_1^2 + |v(s)|_0^2 C(T)\big( 1+ \sup_{ 0\leq r \leq T}
|u(r)|^{4q}_\infty\big)\\
\leq&|v(s)|_1^2 + |v(s)|_1^2 C(T)\big( 1+ \sup_{0\leq r \leq T}
|u(r)|^{4q}_\infty\big)\ .
\end{align*}
This in turn produces
\begin{equation}\label{eq:KVnormRD}
 \E   \sup_{0\leq s, t \leq T}|J_{s,t}|_{\V \rightarrow \V}^{2p} \leq C(T,p)\Big(1
 + \E \sup_{ 0\leq r \leq T}
|u(r)|^{4pq}_\infty\Big)  < \infty,
\end{equation}
for any $p > 0$ and $u_0 \in \H_0$. Since all of the operators on the right-hand side of
the governing equation are self-adjoint in this example, the estimates analogous to
\eqref{eq:JHnormRD} and \eqref{eq:KVnormRD} hold for $K_{s,t}$.

Using the estimates used to produce \eqref{eq:JHnormRD}, it is
straightforward to see that  there exists a $K(T)>0$ and $C(T)>0$ such that
for all $r <s \leq T$,
\begin{equation*}
   |K_{s,T} -K_{r,T}|_{\H \rightarrow\H} \leq \big(e^{K(s-r)} -1\big) e^{K T}
   \leq C |s-r|,
\end{equation*}
which proves the estimate \eqref{eq:assumption_on_Lip_of_K}.

\subsubsection{Verification of Assumption \ref{a:J2a}}
Equation \eqref{eq:JmomentWSimple} has already been verified above
since $g_k \in \V^1$. To see the second estimate, observe that for
$\phi \in \H$ 
\begin{align*}
  J_{s,t} \phi = e^{\Delta (t-s)}\phi  + \int_s^t e^{\Delta (t-r)}
    DN(u_r)J_{r,t}\phi dr,
\end{align*}
and hence we have that for $0\leq s< t\leq T$,
\begin{align*}
  | J_{s,t} \phi |_1 &\leq |e^{\Delta(t-s)}|_{\H\rightarrow \V^1} |\phi|_0
  + \int_s^t  |e^{\Delta(t-r)}|_{\H\rightarrow \V^1} |DN(u_r)|_{\H
    \rightarrow \H} |\phi|_0\\
  &\leq  C \Big(\frac{1}{\sqrt{t-s}} + (1+ \sup_{0\leq r\leq
    T}|u_r|_\infty) \sqrt{t-s} \Big) |\phi|_0\ .
\end{align*}
When Combined with \eqref{eq:infBoundRD}, we obtain that for every $p \geq 1$ and $u_0 \in \H_0$
\begin{equation*}
  \E  \sup_{0\leq s < t < T} (t-s)^\frac12| J_{s,t}|_{\H \rightarrow \V^1}^p \leq C(p,T,u_0)\ .
\end{equation*}

In light of the preceding calculations, we have proven the following result.
\begin{lemma} In the above setting, Assumptions \ref{a:J1}, \ref{a:K1},
  \ref{a:D2}, \ref{a:solReg}, \ref{a:K2a}, and \ref{a:J2a} hold with
  $\HH=\H$, $\HH_0=\H_0$ and $\VV=\V^1$ and $N \in \Poly_1(\V^1,\H)$. Hence, the
  conclusions of Theorem \ref{th:smoothnessA} hold.
\end{lemma}

\subsection{2D Navier Stokes Equation}

Consider the vorticity formulation of the
Navier--Stokes equation in 2D given by:
\begin{equation}\label{eq:NS}
\left\{
  \begin{aligned}
  d w &= \nu \Delta w \,dt + B(\mathcal{K} w,w)\,dt + f(t)\,dt +
  \sum_{j=1}^d g_k dW_k(t)\\
  w(t)&=w_0 \in \H=L^2_0\big([0,2\pi]^2\big),
  \end{aligned}
  \right.
\end{equation}
where $B(u,v) = -(u\cdot \nabla)v$ is the usual Navier--Stokes
nonlinearity, and $\mathcal{K}$ is the Biot--Savart integral operator
which is defined by $u=\mathcal{K} w$ when $w =\nabla \wedge u$ (see
\cite{b:MajdaBertozzi02,MattinglyPardoux05} for more details). We
denote by $L^2_0$ the Hilbert space of square-integrable functions on
$[0,2\pi]^2$ which are periodic and have spatial mean zero. As before,
we form the space $\V^s$, $s \in \RR$,  from $\H=L^2_0$ and
$L=(-\Delta)$.  We assume that $f(t)$ is a bounded function in $\V^1$,
$g_k \in \V^2$.

\begin{lemma}
  In the above setting, Assumptions \ref{a:J1}, \ref{a:K1},
  \ref{a:D2}, \ref{a:solReg}, \ref{a:K2a}, and \ref{a:J2a} hold with
  $\HH=\H$, $\HH_0=\H_0$, and $\VV=\V^1$. Additionally, the map $u
  \mapsto B(\mathcal{K}u,u) \in \Poly_1(\V^1,\H)$. Hence, the
  conclusions of Theorem \ref{th:absolute_continuityA} and
  \ref{th:smoothnessA} hold for equation \eqref{eq:NS}.
\end{lemma}
\begin{proof}
  We begin by proving that $B(\mathcal{K}u,u) \in
  \Poly_1(\V^1,\H)$. To do so we use the basic facts that $|B(u,v)|_0
  \leq C|u|_1| |v|_1$ and that $|\mathcal{K} u|_1 = |u|_0$ (see for
  instance \cite{b:CoFo88}). Then
  \begin{align*}
    |B(\mathcal{K} u,v)|_0 \leq C |u|_0 |v|_1,
  \end{align*}
  which proves the first result. Assumptions \ref{a:J1} and
  \ref{a:K1} then follow from Corollary \ref{c:A1and2OK} or from Proposition
  2.1, Proposition 2.2 of \cite{MattinglyPardoux05}. The
  existence of solutions to \eqref{eq:NS} can also be found in
  \cite{b:Fl94}. Assumptions \ref{a:D2}, \ref{a:solReg} and
  \ref{a:K2a} follow from Corollary A.2 and Lemma B.1 of
  \cite{MattinglyPardoux05}. The fact that $w(t) \in
  \D^\infty_{\V^1}$ (see Section \ref{secMalliavin} for the
  definition) is also proved in Lemma C.1 of
  \cite{MattinglyPardoux05}.
\end{proof}
Lastly we we give a fairly weak condition ensuring that the system is
formally H\"ormander. The following result is a direct consequence of
Corollary 4.5 from \cite{b:HairerMattingly04}.

\begin{lemma}\label{l:nsGinf}
  Let $\mathcal{Z}_0$ be a subset of $\Z^2/(0,0)$ such that the
  following conditions hold:
  \begin{enumerate}
  \item Integer linear combinations of $\mathcal{Z}_0\cap
    (-\mathcal{Z}_0) $ generate $\Z^2$.
  \item There exist two elements of $\mathcal{Z}_0$ with non-equal
    Euclidean norm.
  \end{enumerate}
Then $\Hbb_\infty=\H(\eqdef\HH)$ if 
\begin{equation*}
  \{ \cos( k\cdot x), \sin( k \cdot x) : k \in  \mathcal{Z}_0 \}
  \subset \Hbb_1\eqdef\spn\{g_1,\cdots,g_d\}\ . 
\end{equation*}
\end{lemma}
\begin{remark}
  This result is very similar to one of the principal results in
  \cite{MattinglyPardoux05}. One difference is that we do not
  require that the set of forcing functions consists of $\sin$ or $\cos$ but
  only that the span of the forcing functions contains the needed
  collection of $\sin$ and $\cos$. For a discussion of what happens
  when the conditions in Lemma \ref{l:nsGinf} fail, see
  \cite{b:HairerMattingly04}.
\end{remark}

\section{Malliavin Calculus}\label{secMalliavin}
Since all of our results use techniques from Malliavin calculus, we
give a quick introduction, mainly to fix notation. For a longer
introduction see \cite{MattinglyPardoux05}, for even more background
see e.g.~\cite{b:Nualart95, b:Bell87}.

First, we define the Malliavin derivative of $u(T)$ in the direction $h\in
L^2([0,T],\R^d)$ as
\begin{equation*}
  \MD(u(T))(h)\eqdef\HH-\lim_{\eps\to0}\frac{\Phi(W+\eps H)-\Phi(W)}{\eps},  
\end{equation*}
where $H(t)=\int_0^th(s)ds$. 
It is easy to verify that under Assumption \ref{a:J1} the derivative
$\MD(u(T))(h)$ is well-defined for any $h \in L^2([0,T],\RR^d)$ and that
\begin{equation*}
\MD(u(T))(h)=\int_0^T J_{s,T}Gh(s)ds.  
\end{equation*}
The Malliavin covariance operator $\MM(u(T)):\HH\to\HH$ is defined by
\begin{equation*}
  \ip{\MM(u(T))\phi}{\phi}\eqdef\sum_{k=1}^{d}\int_0^T\ip{J_{s,T}g_k}{\phi}^2ds.  
\end{equation*}
(We shall often write $\MM=\MM(u(T))$ for brevity). It is clearly
nonnegatively definite.  Its finite-dimensional projection on the
space $S$ is given by the Malliavin matrix
\begin{equation}
  \MM_{ij}^\Pi\eqdef\MM_{ij}(\Pi u(T))=\sum_{k=1}^d \int_0^T\ip{J_{s,T}g_k}{\psi_i}
  \ip{J_{s,T}g_k}{\psi_j}ds,\quad i,j=1,\ldots,N, 
  \label{eq:Malliavin_operator}
\end{equation}
where $\psi_1,\ldots,\psi_N$ is an orthonormal basis in $S$.

Notice that the definition in \eqref{eq:Malliavin_operator} involves
solving a continuum of linear systems (one for each $s\in[0,T]$). It is
more convenient to work with the following representation
\begin{equation}
  \label{eq:Malliavin_operator_alternative_representation}
  \ip{\MM(u(T))\phi}{\phi}=\sum_{k=1}^{d}\int_0^T\ip{g_k}{K_{s,T}\phi}^2ds
\end{equation}
which involves solving only one linear system. This representation 
follows from the relation $K_{T,T}\phi=\phi$ and the next lemma:

\begin{lemma} Assume that Assumptions \ref{a:J1} and \ref{a:K1} hold.
  Then for any $\phi,\psi \in\HH$ and $0 < s <t \leq T$, the map
  \begin{equation*}
     r\mapsto\ip{J_{s,r}\phi}{K_{r,t}\psi }
  \end{equation*}
  from $[s,t]$ into $\R$ is constant.
\end{lemma}
\begin{proof} The following essentially recapitulates the proof of
  Proposition 2.3 from \cite{MattinglyPardoux05}.  Set
  $v(r)=J_{s,r}\phi$ and $w(r)=K_{r,t}\psi $. Since $v,w\in
  L^2([s,t],\VV)$ and their time derivatives $v',w'\in
  L^2([s,t],\VV')$, we may apply integration by parts (see Theorem 2
  from~\cite[p.477]{b:DautrayLions88}):
  \begin{multline*}
    \ip{v(r_1)}{w(r_1)}-\ip{v(r_0)}{w(r_0)}=\int_{r_0}^{r_1}
    \left[\ip{v'(r)}{w(r)}+\ip{v(r)}{w'(r)}\right] dr
    \\=\int_{r_0}^{r_1}
    \left[\ip{(DF)(u(r))J_{s,r}\phi}{K_{r,T}\psi
      }-\ip{J_{s,r}\phi}{(DF)^*(u(r))K_{r,T}\psi }\right] dr=0,
  \end{multline*}
  for all $r_0<r_1$.
\end{proof}

\subsection{Higher Malliavin Derivatives} 
The existence of a smooth density requires control of higher Malliavin
derivatives, which we now introduce. For $n\in\N$, $s_1,\cdots,s_n \in
[0,T]$, and
$h_1,\ldots,h_n\in\R^d$, we define
\begin{equation}\label{eq:higherMalDer}
  \MD^{(n)}_{s_1,\ldots,s_n}(u(t))(h_1,\ldots,h_n)\eqdef J^{(n)}_{s_1,\ldots,s_n;t}(Gh_1,\ldots,Gh_n),  
\end{equation}
where $J^{(n)}_{s_1,\ldots,s_n;t}(\phi_1,\ldots,\phi_n)$ is the solution of the $n$-th
equation in variations defined below.

The first variation of equation \eqref{eq:spde} is 
\begin{equation*}
\left\{
  \begin{aligned}
\frac{\partial}{\partial t}J^{(1)}_{s;t}\phi=&DF(u(t)) J^{(1)}_{s;t}\phi,\quad t>s,\\
J^{(1)}_{s;t}\phi=&\phi,\quad t\leq s 
  \end{aligned}\right.
\end{equation*}
for all $\phi\in\VV$. Obviously, $J^{(1)}_{s;t}\phi=J_{s,t}\phi$, where the latter
is introduced in \eqref{eq:J}.

To write down the equations for the higher order variations, we need
some additional notation.  Suppose we have vectors $(s_1,\ldots,s_n)\in\R^n$ and
$(\phi_1,\ldots,\phi_n)\in\VV^n$.  For a subset $I=\{n_1<\ldots<n_{|I|}\}$ of $\{1,\ldots,n\}$
(here $|I|$ means the number of elements in $I$) we denote
$s(I)=(s_{n_1},\ldots,s_{n_{|I|}})$ and $\phi(I)=(\phi_{n_1},\ldots,\phi_{n_{|I|}})$.

Now for $n\geq2$, $s_1,\cdots,s_n$, and  $\phi_1,\ldots,\phi_n\in\VV$, the $n$-th equation in variations
is given by

\begin{align}\label{eq:higher_variations}
  \frac{\partial}{\partial t}J^{(n)}_{s_1,\ldots,s_n;t}(\phi_1,\ldots,\phi_n)=&DF(u(t))
  J^{(n)}_{s_1,\ldots,s_n;t}(\phi_1,\ldots,\phi_n)\notag\\
  &\quad+G^{(n)}_{s_1,\ldots,s_{n};t}(u(t))(\phi_1,\ldots,\phi_n),\quad t>\lor s,\\
  J^{(n)}_{s_1,\ldots, s_n; t}(\phi_1,\ldots,\phi_{n})=&0, \quad t\leq \lor s,\notag
\end{align}
where $\lor s =s_1\lor\ldots\lor s_n$, and for $n\in\N$,
\begin{multline}
  \label{eq:G^{(n)}}
  G^{(n)}_{s_1,\ldots,s_{n};t}(u(t))(\phi_1,\ldots,\phi_n)\\= \sum_{\nu=1}^n\sum_{I_1,\ldots,I_\nu}
  D^{(\nu)}F(u(t))\left(J_{s(I_1);t}^{|I_1|}\phi(I_1),\ldots,J_{s(I_\nu);t}^{|I_\nu|}\phi(I_\nu)\right)
  \\=\sum_{\nu=2}^{m\land n}\sum_{I_1,\ldots,I_\nu}
  D^{(\nu)}N(u(t))\left(J_{s(I_1);t}^{|I_1|}\phi(I_1),\ldots,J_{s(I_\nu);t}^{|I_\nu|}\phi(I_\nu)\right).
\end{multline}
Here $m$ is the degree of the polynomial $F$, and the inner sum is
taken over all partitions of $\{1,\ldots,n\}$ into disjoint nonempty
sets $I_1,\ldots,I_\nu$ (we do not distinguish two partitions obtained
from each other by a permutation).  The upper limit in the outer sum
can be changed to $m\land n$ since the derivatives of $F$ of order
higher than $m$ vanish. The lower limit can be changed to $2$ since
there are no admissible partitions for $\nu=1$. Since  $N$ has all of
the non-linear terms in the equation we can replace $F$ with $N$.

Variation of constants for~\eqref{eq:higher_variations} gives
\begin{equation}
  \label{eq:integral_for_higher_J_1}
  J^{(n)}_{s_1,\ldots,s_n;t}(\phi_1,\ldots,\phi_n)
  =\int_{\lor s}^t J_{r,t}G^{(n)}_{s_1,\ldots,s_{n};r}(u(r))(\phi_1,\ldots,\phi_n)dr,
\end{equation}
for $n\geq 2$.

We say that $u(t)\in\D^\infty_Y$ for some Banach space $Y$ if for all
$n\in\N$, and all $h_1,\ldots,h_{n}\in\R^d$,
\begin{equation}
\label{eq:def_of_D_infty}
\E |\MD^{(n)}_{s_1,\ldots,s_n}(u(t))(h_1,\ldots,h_{n})|^p_Y<\infty,\ \mbox{for all}\
p>1.
\end{equation}

\begin{lemma} \label{l:UMalSmooth}Under assumptions \ref{a:solReg} and
  \ref{a:J2a}, for all $n\in\N$,
\begin{equation*}
  \E\;\sup_{\phi_k \in \{g_1,\cdots,g_d\}}\sup_{s,r}
  \|J^{(n)}_{s_1,\ldots,s_n;r}(\phi_1,\ldots,\phi_n)\|^p<\infty
  \quad \mbox{for all\ } p\geq1,    
\end{equation*}
and hence $u(T)$ belongs to $\D^\infty_\VV$.
\end{lemma}

\begin{proof}
  The fact that $u(T)$ belongs to $\D^\infty_\VV$ follows immediately
  from the first part of the Lemma when combined with
  \eqref{eq:higherMalDer}. The first claim will follow by induction.

For $n=1$ the statement for $J^{(1)}$ follows directly from
\eqref{eq:JmomentWSimple} in
Assumption~\ref{a:J2a}.
Let us fix $n\geq2$ and suppose that the statement holds true for all
positive integers less than $n$. Take any $0<s_1,\ldots,s_n<r<T$ such
that $\lor s =\TTT$. In the interest of notational compactness we write
\begin{equation*}
  \sum_{\nu,I,j} \quad \text{for}\quad \sum_{\nu=2}^{m\land
    n}\sum_{I_1,\ldots,I_\nu}\sum_{j=\nu}^m\;.
\end{equation*}
  Then by~\eqref{eq:integral_for_higher_J_1}, we have
\begin{multline*}
  \|J^{(n)}_{s_1,\ldots,s_n;r}(\phi_1,\ldots,\phi_n)\|\leq
  \int_{\lor s}^T
  \left\|J_{r,T}G^{(n)}_{s_1,\ldots,s_{n};r}(u(r))(\phi_1,\ldots,\phi_n)\right\| dr\\
  \leq \int_{\lor s}^T \sum_{\nu,I,j}  \left\|J_{r,T}
    D^{(\nu)}N_j(u(r))\left(J_{s(I_1);r}^{|I_1|}\phi(I_1),\ldots,
      J_{s(I_\nu);r}^{|I_\nu|}\phi(I_\nu) \right)\right\| dr \\\leq
  m!  \sum_{\nu,I,j} \int_{\lor s}^T \left\|J_{r,T}
    N_j\left(u(r)^{\otimes 
        j-\nu},J_{s(I_1);r}^{|I_1|}\phi(I_1),\ldots,J_{s(I_\nu);r}^{|I_\nu|}
      \phi(I_\nu)\right)\right\| dr \\\leq m!
  \sum_{\nu,I,j}\left[\sup_r(T-r)^\alpha|J_{r,T}|_{\HH\to\VV}\right]
  \sup_{r,s} \|u(r)\|^{j-\nu}
   \prod_{l=1}^\nu \|J_{s(I_l);r}^{|I_l|}\phi(I_l)\|
  \int_{\lor s}^T\frac{1}{(T-r)^{\alpha}}dr.
\end{multline*} Since the r.h.s. is bounded uniformly in $s$, $r$, and
the choice of $\phi_k$, we can take the supremum over all of them. Next,
taking the expectation of both sides, we use H\"older's inequality to
split the products. The estimates \eqref{eq:JmomentWSimple},
\eqref{eq:pth_momment_JSmooth}, and
\eqref{eq:assumption_moment}, and the induction hypothesis imply that
all the moments of the r.h.s. are finite, and we are done.
\end{proof}

\section{ General Results}
\label{sec:GeneralResults}
We now give the proof of the main results of this article. They are
generalizations of the results given in Sections \ref{sec:ExtDenBasic}
and \ref{sec:smoothDenBasic}. All of our examples fit into the
framework of the previous sections. However, for completeness and to
emphasize the connection with the standard finite-dimensional results,
we will prove the more general results in this section, which imply the
results previously stated.

\subsection{Existence of densities}\label{sec:ExistenceOfDensities}

To understand how the randomness spreads through the phase space, we
now introduce an increasing collection of sets which characterize the
directions excited.

The Lie bracket of two Fr\'echet-differentiable vector fields
$A,B:\VV\to\VV'$ is a new vector field
\begin{equation*}
  [A,B](x) \eqdef (DA)(x)B(x) - (DB)(x)A(x)\in\VV',  
\end{equation*}
defined for all $x\in\VV$ when it makes sense (i.e. when
$A(x),B(x)\in\VV$). In the interest of notational brevity, we will write
\begin{equation}
  \label{eq:R}
  \bl A_1,A_2,\ldots,A_n\br
  (x)\eqdef[\ldots[[A_1,A_2],\ldots],A_n](x)\ .
\end{equation}

Next, we define the set $\Ab$ of admissible vector fields which will
play an essential role in the forthcoming iteration scheme. To do so
we fix a time $\Tt \in [0,T)$ and recall the process $X(t)=u(t)-GW(t)$
defined earlier. Notice that $X(t)$ can also be written as
\begin{equation}\label{eq:X}
  X(t,\omega)=X(t)=u(0)+\int_0^tF(u(s))ds+\int_0^tf(s)ds=u(t)-GW(t),
\end{equation}
and hence $X(t) \in C\big( [0,T],\HH\big)\cap C\big( (0,T], \VV\big)$ almost surely.

\begin{definition}
  $\Ab$ is  the set of all
  polynomial vector fields $Q:\VV\to\VV'$ such that with probability one the
  following conditions hold:
  \begin{enumerate}
  \item $Q(X(t))\in L^2\big([\Tt,T],\VV\big)$,
  \item $\frac{d}{dt}Q(X(t))\in L^2\big([\Tt,T],\VV'\big)$,
  \item $[F,Q]$ is a continuous polynomial from $\VV \to\HH$.
  \end{enumerate}
\end{definition}

For any $\omega\in\Omega$ and any positive integer $n$, we introduce a
set $\Hb_n$ of smooth vector fields $Q:\VV\to\VV'$.  For $n=1$, we set
$\Hb_1=\spn\{g_1,\ldots,g_d\}$.  For $n>1$, $\Hb_n$ is defined
recursively from $\Hb_{n-1}$:
\begin{equation}\label{a:HaDef}
  \Hb_n\eqdef\spn\biggl(\Hb_{n-1}\cup \bigcup_{Q\in \Hb_{n-1} \cap \Ab}
    \;\bigcup_{i=0}^\infty\; \bigcup_{k_1,\ldots,k_i} \{ \bl
    F,Q,g_{k_1}, \ldots ,g_{k_i} \br \} \biggr). 
\end{equation}

Now we introduce $\Hb_\infty=\spn\big(\bigcup_n\Hb_n\big)$ and for
$n\in\N\cup\{\infty\}$ define

\begin{equation*}
  \Hb_n(x)\eqdef\{Q(x):\ Q\in\Hb_n\}.  
\end{equation*}

\begin{theorem}\label{th:absolute_continuityB} Assume that Assumptions
  \ref{a:J1}, \ref{a:K1}, and \ref{a:D2} hold.  Suppose that $S$ is a
  finite-dimensional linear subspace in $\HH$.  If, in addition, $S$ is
  a subspace of $\Hb_\infty(X(T))$ with probability 1, then the
  distribution of the orthogonal projection $\Pi_S u(T)$ on $S$ is
  absolutely continuous with respect to the Lebesgue measure on~$S$.
\end{theorem}

We will see in remark \ref{rm:relaxD2} that the above theorem holds
under a slightly relaxed version of Assumption \ref{a:D2}.
The following lemma shows that Theorem \ref{th:absolute_continuityA}
is implied by Theorem \ref{th:absolute_continuityB} given above. Its
proof will be given after the proof of 
Theorem~\ref{th:absolute_continuityB}.
\begin{lemma}
  \label{l:AtoBabs}
  Under Assumptions \ref{a:J1} and \ref{a:K1}, $\Hbb_n \subset \Hb_n$
  for all $n$.
\end{lemma}

\begin{proof}
  We shall proceed by induction.
  First, notice that $\Hbb_1=\Hb_1$ and that all of the elements of
  $\Hbb_1$ are constant. Now our induction hypothesis
  will be that for some $n$ we have $\Hbb_{n-1}\subset\Hb_{n-1}$ and
  that all vector fields in $\Hbb_{n-1}$ are constant.
  
  It is sufficient to show that if $h=N_m(g,g_{k_1},\ldots,g_{k_{m-1}})\in\VV$ for
  some $g\in\Hbb_{n-1}\cap\VV\cap\Dom(L)$ and $k_i\in\{1,\ldots,d\}$, then 
  $h$ is constant in $\VV$ (which is trivial), and there is a $Q\in\Hb_{n-1}\cap\Ab$ such that
  $h=[F,Q,g_{k_1},\ldots,g_{k_{m-1}}]$.
  
  To prove the latter, we can choose $Q=g/m!$. Then
  Lemmas
  \ref{lm:brackets_with_constant_fields} and
  \ref{lm:higher_derivatives_of_multilinear} imply: 
  \begin{equation*}
    [F,Q,g_{k_1},\ldots,g_{k_{m-1}}]=N_m(g,g_{k_1},\ldots,g_{k_{m-1}})=h.
  \end{equation*}
  We shall now check that $Q$ or, equivalently, $g$ belongs to $\Hb_{n-1}\cap\Ab$.
  First, notice that $g\in\Hb_{n-1}$ by the induction
  hypothesis. Next, $g\in\Ab$ since i) $g\in\VV$, ii) $\frac{d}{dt}g=0$, 
  and iii) Lemma
  \ref{lm:brackets_with_constant_fields} shows that
  $[F,g]=(DF)(x)(g)=F(g\otimes x^{\otimes m-1})=Lg+N(g\otimes
  x^{\otimes m-1})$, which is continuous from $\VV \to \HH$
  by the assumptions on $N$, since $Lg$ is a constant in $\HH$ due to
  $g\in\Dom(L)$. 
\end{proof}

\subsection{Smoothness of densities}
\label{ssec:smoothDensities}
We now introduce a second sequence of sets $\Ha_n$
of vector fields from $\VV$ to $\HH$. The  $\Ha_n$ play the analogous
role in our smoothness of density result as the  $\Hb_n$ played in the
existence of density result.  We begin by defining a slightly modified
version of the set of admissible vector fields $\A$  used in the last
section. Let
\begin{equation}\label{defAa}
  \Aa\eqdef\Big\{ Q:\VV \to \HH : Q \in \Ab \cap {\Poly}_1(\VV,\HH) \Big\}.
 \end{equation}

Given a collection of functions $\mathcal{C}$ we define the symmetric
convex hull, denoted $\sch(\mathcal{C})$, by
\begin{equation*}
  \sch (\mathcal{C}) \eqdef \biggl\{ \sum \alpha_i f_i : f_i \in
    \mathcal{C}, \alpha_i \in \RR, \text{ and } \sum_i |\alpha_i| \leq 1 \biggr\}. 
\end{equation*}

For $n=1$, we set $\Ha_1=\sch(g_1,\ldots,g_d) \subset \Hb_1$.  For
$n>1$, we construct $\Ha_n$ from~$\Ha_{n-1}$.

We set
\begin{equation}\label{eq:HaDef}
  \Ha_{n}\eqdef\sch\bigg(\Ha_{n-1}\cup\bigcup_{Q\in\Ha_{n-1}\cap\Aa}
    \;\bigcup_{i=0}^\infty\;\bigcup_{k_1,\ldots,k_i} \{\bl
    F,Q,g_{k_1},\ldots,g_{k_i}\br\} \bigg)\ .     
\end{equation}

\begin{theorem}\label{th:smoothness}
  Assume that Assumptions \ref{a:solReg}, \ref{a:K2a} and \ref{a:J2a}
  hold.  Let $S$ be a deterministic finite-dimensional subspace of
  $\VV$ such that for some $n$ and some $\delta>0$
\begin{equation}
  \label{eq:approximation}
 \Lambda_p^*(u_0,T) \eqdef \E\left[\inf_{\substack{\|\phi\|\le 1\\
 \|\Pi_S
        \phi\|\ge \delta}}\sup_{Q\in \Ha_n(X(T))}\Bigl|\langle
    \phi,Q(X(T))\rangle\Bigr|\right]^{-p}<\infty, 
\end{equation}
for all $p\geq1$.
Then the density of $\Pi_S u(T)$ with respect to Lebesgue measure (whose
existence is guaranteed by Theorem~\ref{th:absolute_continuityB}) is a
$C^\infty$-function on $S$.
\end{theorem}

The next lemma shows that Theorem \ref{th:smoothnessA} follows from
Theorem \ref{th:smoothness}.
\begin{lemma} Recall the definition of $\Hbb_n$ from \eqref{def:G}.
  If $S \subset \Hbb_n$ then the condition in \eqref{eq:approximation}
  holds for this $n$. In fact, there exists a subset of constant vector
  fields $\Ha_n^\prime\subset \Ha_n$ such that
  \begin{equation*}
    \inf_{\substack{\|\phi\|\le 1 \\ \|\Pi_S
        \phi\|\ge \delta}}\sup_{Q\in \Ha_n^\prime}\Bigl|\langle
    \phi,Q\rangle\Bigr| >0
  \end{equation*}
for some $\delta >0$.
\end{lemma}

\section{Proof of General Results}\label{sc:genResults}

\subsection{Absolute continuity}\label{sc:absolute_continuityB}

Theorem~\ref{th:absolute_continuityB} will be implied by the following
standard result from Malliavin calculus (see~\cite[p.86, Section
2.1]{b:Nualart95}; it is straightforward to check that the definitions of
the Malliavin derivative and matrix given in~\cite{b:Nualart95} are
equivalent to ours):
\begin{theorem}\label{thm:absolute-continuityNualart}
  Suppose the following conditions are satisfied for a
  finite-dimensional random vector $Y$:
  \begin{enumerate}
  \item $\E|\MD(Y)(h)|^2\leq\infty$, for all $h\in L^2\big([0,T],\R^d\big)$.\label{i:abs1}
  \item The Malliavin matrix $\MM(Y)$
    is invertible a.s.\label{i:abs2}
  \end{enumerate}
  Then the law of $Y$ is absolutely continuous with 
  respect to the Lebesgue measure.  \label{th:density_nualart}
\end{theorem}

\begin{proof}[Proof of Theorem~\ref{th:absolute_continuityB}]
  Condition~\ref{i:abs1}) of Theorem~\ref{th:density_nualart} follows
  from~\eqref{eq:assumption_on_2nd_moment_of_J}. 
  
  To verify condition~\ref{i:abs2}), it is sufficient to prove that
  \begin{equation}
    \label{eq:sufficient_for_inverse}
    \PP\Big\{\Ker\MM\cap \Hb_{\infty}(X(T)) \neq \{0\} \Big\}=0,
  \end{equation}
  where
  \begin{equation*}
    \Ker\MM =\Big\{\phi\in\VV:\ip{\MM\phi}{\phi}=0\Big\}.
  \end{equation*}
  This property is implied by
  \begin{equation*}
    \PP\Big\{\Ker\MM\perp \Hb_\infty(X(T))\Big\}=1    
  \end{equation*}
  or, equivalently, by
  \begin{equation*}
    \PP\Big\{\Ker\MM\perp \Hb_n(X(T))\ \mbox{for all}\ n\in\N\Big\}=1,
  \end{equation*}
  which in turn follows from
  \begin{equation}
    \PP\Big\{\Ker\MM\perp Q(X(T)),\ \mbox{for all}\ Q\in\Hb_n, n\in\N\Big\}=1,
    \label{eq:Ker_orthogonality_G}
  \end{equation}
  and the fact that $\Hb_\infty(X(T))$ is generated by $Q(X(T)),Q\in\Hb_n, n\in\N$. 
  Relation~\eqref{eq:Ker_orthogonality_G} is a consequence of the
  following statement which we will prove below. There is a set $\Omega'$ with $\PP(\Omega')=1$
  such that for all $\omega\in\Omega'$, all $\phi\in\Ker\MM$,  every
  $n\in\N$, each  $Q\in\Hb_n$ and all
  $s\in[\Tt,T]$, we have that
  \begin{equation}
    \ip{Q(X(s))}{K_{s,T}\phi}=0
    \label{eq:Ker_orthogonality_G_all_s}
  \end{equation}
  where $\Tt$ was the time fixed at the start of Section \ref{sec:ExistenceOfDensities}.

  This statement will be proved by induction in $n$.  For $n=1$ it
  follows directly from the
  representation in~\eqref{eq:Malliavin_operator_alternative_representation}.
  The induction step is provided by the next lemma, whose proof will
  complete the proof of the present result.
\end{proof}

\begin{lemma}\label{lm:induction_step} There exists a set $\Omega'$ of probability~1
  such that for all $\omega$ in this set $\Omega'$, the following
  implication holds true:

Let $Q:\VV\to\VV'$ be a polynomial vector field
in $\Ab$. Then  for any $t_0 \in [\Tt,T]$,
\begin{equation}
  \langle Q(X(s)),K_{s,T}\phi\rangle=0,\quad s\in[t_0,T],
\label{eq:orthogonality}
\end{equation}
implies that 
\begin{equation*}
  \left\langle[F,Q,g_{k_1},g_{k_2},\ldots,g_{k_i}](X(s)), K_{s,t}\phi
  \right\rangle=0, \quad s \in [t_0,T],  
\end{equation*}
for any $i\geq 0$ and $k_j \in \{1,\ldots,d\}$. 
\end{lemma}
\begin{proof}
  By Theorem 2 from~\cite[p.477]{b:DautrayLions88}, we can
  differentiate~\eqref{eq:orthogonality} with respect to $s$.
  Equation~\eqref{eq:adjoint_equation} implies
\begin{multline*}
  0=\frac{\partial}{\partial s} \langle Q(X(s)),K_{s,T}\phi\rangle\\ = \langle
  (DQ)(X(s))F(u(s)),K_{s,T}\phi\rangle-\langle
  Q(X(s)),(DF)^*(u(s))K_{s,T}\phi\rangle\\
  =\langle (DQ)(X(s))F(u(s))- DF(u(s))Q(X(s)),K_{s,T}\phi\rangle.
\end{multline*}
Fix $s$ and $X(s)$ and notice that the vector field
\begin{equation*}
  R(y)=DQ(X(s))F(y)- DF(y)Q(X(s))
\end{equation*}
is well-defined and a polynomial from $\VV \to \HH$. We also have
\begin{equation*}
  [R,g_{k_1},g_{k_2},\ldots,g_{k_i}](X(s))=
  -[F,Q,g_{k_1},g_{k_2},\ldots,g_{k_i}](X(s)).
\end{equation*}
Hence by Lemma~\ref{lm:expanding_polylynear}, for $s \in [t_0,T]$ and
some $n$
\begin{align}\label{eq:eqZeroLine}
  0=&\langle R(X(s)),K_{s,T}\phi\rangle \\=&-
  \sum_{i=0}^{n}\sum_{k_1,\ldots,k_i}\langle
  [F,Q,g_{k_1},\ldots,g_{k_i}](X(s)),K_{s,T}\phi\rangle W_{k_1}\ldots
  W_{k_i}.\notag
\end{align}
Observe that each of the inner products is a continuous function of
time. This follows from the almost sure continuity in $\HH$ of the two
arguments of the inner products.  The brackets, by virtue of being in
$\Ab$, are continuous from $\VV \to \HH$, and
$X(s)$ is continuous in $\VV$ on $[\Tt,T]$ almost surely by assumption. Hence, if
$Y(s)= [F,Q,,g_{k_1},\ldots,g_{k_i}](X(s))$, then $Y(t)$ is in
$C([\Tt,T],\HH)$. By assumption, we know that $K_{s,T}\phi$ is
$C([\Tt,T],\HH)$ almost surely. For $\Tt \leq s<r \leq T$ we have 
\begin{align*}
  |\langle Y(r),K_{r,T}\phi\rangle - \langle Y(s),K_{s,T}\phi\rangle| &\leq
  |\langle Y(r)-Y(s),K_{s,T}\phi\rangle| +   |\langle
  Y(s),K_{r,T}-K_{s,T}\phi\rangle|\\
  &\leq  |Y(r)-Y(s)| |K_{s,T}\phi| + |Y(s)|
  |K_{r,T}-K_{s,T}|,
\end{align*}
and thus conclude that $\langle Y(r),K_{r,T}\phi\rangle$ is continuous in $r$.
The proof of the result is now completed using Theorem
\ref{th:quadratic_variation_of_polynomial}. 
\end{proof}

\subsection{Smoothness of the density}\label{sc:smoothness}

Theorem~\ref{th:smoothness} will follow from the following
classical result from Malliavin Calculus (see for example
\cite[Corollary 2.1.2]{b:Nualart95}) which is a strengthening of
Theorem \ref{th:density_nualart} which was used to prove the existence
of a density:
\begin{theorem}\label{th:smoothness_nualart}
  Suppose that $\Pi$ is the orthogonal projection onto some 
  finite-dimensional subspace of $\YY$ and the following conditions hold:
  \begin{enumerate}
  \item $\Pi u(T)$ belongs to $\D^\infty_\YY$.
  \item The projected Malliavin matrix $\MM^\Pi=\MM(\Pi u(T))=(\MM_{ij})$ (defined in \eqref{eq:Malliavin_operator} )  satisfies 
    \begin{equation*}
      \E|\det\MM^\Pi|^{-p}<\infty\ \mbox{for all}\ p>1.
    \end{equation*}
  \end{enumerate}
  Then the density of $\Pi u(T)$ with respect to Lebesgue measure on
  $\YY$ exists and is $C^{\infty}$-smooth.
\end{theorem}

We have to check both conditions of this theorem to prove 
Theorem~\ref{th:smoothness}. The first condition
is implied by Lemma~\ref{l:UMalSmooth}, and the second one
follows from the theorem below. For $n \in \N \cup
\{\infty\}$, we define $\Sa_n = \spn( \Ha_n)$.  Here
$\Ha_\infty=\cup_n^\infty \Ha_n$.


\begin{theorem}
  \label{lm:moments_matrix_inverse} Let $\Pi$ be the orthogonal
  projection onto a finite-dimensional subspace of $\Sa_n$ for some $n$. 
  Fix a number $\delta>0$. Let $U=U_\delta=\{\phi\in \VV: \|\phi\|\le1,\ \|\Pi\phi\|\geq \delta \}$. Then for any $p \geq
  1$, there is
  $\eps_0 = \eps_0(p)$ such that
  \begin{equation*}
    \PP\left\{\inf_{\phi\in
        U}\ip{\MM(u(T))\phi}{\phi}<\eps\right\}\leq \eps^p   
  \end{equation*}
  if $\eps\leq\eps_0$.
\end{theorem}

\begin{remark}\label{rm:relaxD2} We notice that Assumption \ref{a:D2}
  can be relaxed.   Specifically, to satisfy the first condition in Theorem
  \ref{thm:absolute-continuityNualart}, we only need second moments of
  the Malliavin derivative of $\Pi u(T)$. We only need Assumption
  \ref{a:D2} to hold with the $|\cdot|$ norm replaced by a norm dual
  to a norm, which is finite on $S$. For instance if $S \subset \VV$
  then \eqref{eq:assumption_on_2nd_moment_of_J} can be replaced by
  \begin{equation*}
       \sup_k \sup_{0 < s<t\leq T}\E|J_{s,t}g_k|_{\VV'}^2\leq J^*(T,u_0).
  \end{equation*}

\end{remark}

\subsubsection{The Proof of Theorem \ref{lm:moments_matrix_inverse} and Associated Results}\label{sec:proofThmFixT0}
The proof of this theorem will use a quantitative version of
Lemma~\ref{lm:induction_step}.  From this point forward, we fix $\TT$
to be the maximum of the two $\TTT$'s given in Assumptions
\ref{a:solReg} and \ref{a:K2a}.

Before stating the result, we need a little notation: For $f:[0,T] \to
\RR$ we define
\begin{equation*}
  \LipT(f)\eqdef \sup_{\TT\leq s< t\leq T}\frac{|f(t)-f(s)|}{t-s}
  \quad\text{and}\quad \supT(f)\eqdef \sup_{\TT\leq  t\leq T} |f(t)|\ .
\end{equation*}
If $f:[0,T] \to \VV$, then by $\LipT|f|$ and $\supT\|f\|$ we mean the
same expressions with the absolute values replaced by the indicated
norm.  When applied to the operator $K_{s,t}$ we mean the same
expressions where $s$ and $t$ vary over all $s,t \in [\TT,T]$ with $s
<t$.  Lastly, we define $\Norm{g}\eqdef\max\big\{\LipT(g),\;
\supT(g)\big\}$,
\begin{equation*}
  \Norm{g}_\HH\eqdef\max\left\{\LipT|g|_{\VV'},\;
    \supT|g|\right\}\quad\text{and}\quad
  \Norm{g}_\VV\eqdef\max\left\{\LipT|g|, \; \supT\|g\|\right\} .  
\end{equation*}

We now give a number of properties of the symmetric convex hull of a
set of functions.
\begin{lemma}\label{l:admisClosed} Recalling the definition of $\Aa$ from equation \eqref{defAa}, let $f_1,\cdots,f_m$ be a collection of polynomial vector fields 
  from $\VV \to \VV'$ with $f_i \in \Aa$ for all $i$. Let
  $\mathcal{C}= \sch(f_1,\cdots,f_m)$. If $ g \in \mathcal{C}$, then
  $g \in \Aa$, and for all $x \in \VV$,
   \begin{equation*}
     \Lip(g)(x) \leq \sup_i \Lip(f_i)(x),
   \end{equation*}
where $\Lip$ is the local Lipschitz constant defined in
$\eqref{eq:lipDef}$ and viewed as a function from $\HH \to \VV'$.
\end{lemma}
\begin{proof}
  Let $\Extr \mathcal{C}$ denote the extreme points of
  $\mathcal{C}$. Clearly, $\Extr \mathcal{C} \subset
  \{f_1,-f_1,\cdots,f_m,-f_m\}$ so it is finite. Being an element of
  $\mathcal{C}$, $g$ is a linear combination of its
  extreme points. Since this set is finite and each $ f_i \in \Aa$,
  we see that $g \in \Aa$. Since $g = \sum \alpha_i f_i$ with $\sum
  |\alpha_i| =1$, we have that
  \begin{equation*}
    \Lip(g) \leq \sum_i |\alpha_i| \Lip(f_i) \leq \sup_i  \Lip(f_i)\;.
  \end{equation*}
\end{proof}
 
\begin{corollary}\label{l:unifLip} For all $n \geq 1$,
   $\Ha_n \subset \Aa$ and $\Ha_n$ is a collection of uniformly locally
   Lipschitz functions from $\VV \to \VV'$, where the $\HH$ norm is used
   on the domain. In particular, there
   is a constant $p(n) \geq 1$  and $C(n) > 0$ so that 
   \begin{equation*}
     \sup_{g \in \Ha_n}\Lip(g)(x) \leq C(1+\|x\|^p), 
   \end{equation*}
for all $x \in \VV$, where $g$ is viewed a polynomial from $\HH
\to \VV'$.
\end{corollary}
\begin{proof}
  Combine Lemma \ref{l:locallipPoly} with Lemma \ref{l:admisClosed}.
\end{proof}

We now give the workhorse lemma which will be used iteratively in the
proof of the main result.

\begin{lemma}\label{lm:induction_step_soft} Recall that $d$ is the
  number of Wiener processes driving the system.  There is a universal,
  positive number $\eps_0(d)$ such that for all $\eps\in(0,\eps_0)$
  there is a set $H(\eps)\subset \Omega$ with the following property:
  If $Q:\VV\to\VV'$ is a vector field in $\Poly_1(\VV,\HH)$
  then for all $\eps\in(0,\eps_0)$ and all $\phi\in\VV$ with
  $\|\phi\|\leq 1$
  \begin{multline*}
    \Bigl\{\ \supT\langle Q(X(s)), K_{s,T}\phi\rangle< \eps, \\
    \max_i\max_{k_1,k_2,\ldots,k_i} \supT \left\langle \bl
      F,Q,g_{k_1},\ldots,g_{k_i}\br(X(t)),K_{s,T}\phi\right\rangle>
    \eps^{8^{-(m+3)}}\Bigr\} \\ \subset
    H(\eps)\cup\Bigl\{\max_i\max_{k_1,k_2,\ldots,k_i}\sup_{\|\phi\|\leq
      1} \Norm{\left\langle \bl F,Q,g_{k_1}, \ldots,g_{k_i}\br(X(s)),
        K_{s,T}\phi \right\rangle}> \eps^{-8^{-(m+3)}}\Bigr\}.
  \end{multline*}
  Here $H(\eps)$ is also universal, depending only on the number
  $d$; $m$ is the degree of the polynomial $F$.

Furthermore, there are universal, positive constants $K_1(d),K_2(d)$, and $\gamma(d)$
such that 
\begin{equation*}
  \PP(H(\eps))\leq K_1e^{-K_2\eps^{\gamma}},
\end{equation*}
for $\eps\in(0,\eps_0(d))$.  
\end{lemma}

With these results stated we return to the proof of Theorem
\ref{lm:moments_matrix_inverse}, postponing the other proofs to the end
of the section.
\begin{proof}[Proof of Theorem \ref{lm:moments_matrix_inverse}]
  First observe that the representation
  \eqref{eq:Malliavin_operator_alternative_representation} implies
  that
  \begin{align*}
    \PP\left\{\inf_{\phi\in U}\ip{\MM(u(T))\phi}{\phi}<\eps\right\}
    &\leq\PP\left\{\inf_{\phi\in U}
      \sum_{k=1}^{d}\int_0^T\ip{g_k}{K_{s,T}\phi}^2ds<\eps\right\}\\ 
    &\leq\PP\left\{\inf_{\phi\in U} 
      \max_{k=1,\ldots,d}\int_{\TT}^T\ip{g_k}{K_{s,T}\phi}^2ds<\eps\right\}
  \end{align*}
  where $\TT$ was again the time fixed at the start of Section \ref{sec:proofThmFixT0}.

  We now need an elementary auxiliary lemma which can be found in
  \cite{MattinglyPardoux05}. We denote by $\Hol_\rho(f)$ 
  the H\"older constant of degree~$\rho$ of a function $f$ (see 
  Section~\ref{sc:more_quantative_estimates} for a precise definition.)
  
  \begin{lemma}\cite[Lemma 7.6]{MattinglyPardoux05} 
    \label{lm:transition_from_l_p_to_l_infty} 
    For any $\eps>0$ and $l>0$, $\int_0^t|f(s)|^lds<\eps$ and
    $\Hol_\rho(f)<c\eps^{-\gamma}$ imply
    $\norm{f}_{L^\infty}<(1+c)\eps^{\frac{\rho-\gamma}{1+l\rho}}$.
  \end{lemma}
  This lemma implies that for a fixed $\phi\in U$ and any
  $l=1,\ldots,d$,
  \begin{multline}
    \label{eq:transition_from_L2_to_Linfty}
    \left\{\max_{k=1,\ldots,d}\int_{\TT}^T\ip{g_k}{K_{s,T}\phi}^2ds
    <\eps\right\}\subset
    \left\{\LipT(\langle g_l,K_{s,T}\phi\rangle)\geq
      \eps^{-1/3}\right\}\\\cup\left\{\supT\ip{g_l}{K_{s,T}\phi}s<\eps^{1/6}\right\},
  \end{multline}
  for  $\eps \in (0,\eps_1]$, where $\eps_1$ is a universal
  constant independent of everything in the problem.
  We also have
  \begin{equation}
    \left\{\LipT(\langle g_l,K_{s,T}\phi\rangle)\geq
      \eps^{-1/3}\right\}\subset\left\{
      |g_l|\; \LipT|K_{s,t}|\geq\eps^{-1/3}\right\}, 
    \label{eq:power_decay1}
  \end{equation}
  where $\LipT|K_{s,t}|=\sup_{|\phi|\le1}\LipT|K_{s,t}\phi|$.
  Notice that the event in the r.h.s. of \eqref{eq:power_decay1} does not depend on $\phi$. Hence
  if we define $g_*~=~\max(1,\sup_i  |g_i|)$, $\Norm{K_{s,t}}_\VV=\sup_{\|\phi\|\le 1}\Norm{K_{s,t}\phi}_{\VV}$, and 
  \begin{align*}
    D_*(R)&=\left\{ g_*\Norm{K_{s,t}}_\VV\geq R\right\},\\
    A_1(\eps)&=\left\{\sup_{\phi\in U}
                \sup_{Q\in \Ha_1} \supT\ip{Q(X(s))}{\phi}<\eps^{1/6}
      \right\},
  \end{align*}
  we have 
  \begin{equation}\label{eq:partialProbBound2}
    \PP\left\{\inf_{\phi\in U}\ip{\MM(u(T))\phi}{\phi}<\eps\right\}
    \leq \PP\Big( D_*(  \eps^{-1/3}) \cup A_1(\eps) \Big).
  \end{equation}
  for all $\eps \in (0,\eps_1]$. Estimates from section 
  \ref{sc:bounds_on_lipschitz_norms} show that
  $D_*(\eps^{-1/3})$ has sufficiently fast decaying 
  probability as $\eps \to 0$, so we need to obtain a good estimate
  on the probability of $A_1(\eps)$. To that end, we define
  \begin{equation*}
    A_i(\eps)=\bigcup_{\phi\in U} A_i(\phi),
  \end{equation*}
  where
  \begin{equation*}
    A_i(\phi)=\left\{\sup_{Q\in \Ha_i\setminus
        \Ha_{i-1}} \supT\ip{Q(X(s))}{\phi}<\eps^{\kappa(i)}\right\},\quad
    i=1,2,\ldots,\quad\eps>0,
  \end{equation*}
  and $\kappa(i)= \frac{1}{6\cdot 8^{(m+3)(i-1)} }$, for
  $i\in\N$. 
  (In this definition, we set $\Ha_0=\emptyset$. Notice that this is
  consistent with the definition of $A_1(\eps)$ given above.)  Next, we define 
  \begin{equation*}
    B_i(\eps)=\bigcup_{\phi\in U} B_i(\phi), 
  \end{equation*}
  where
  \begin{equation*}
    B_i(\phi)=A_{i-1}(\phi)\setminus A_{i}(\phi),\quad i=2,3,\ldots .
  \end{equation*}
  Notice that $A_1= ( A_1 \cap A_2) \cup B_2$, $A_2 = ( A_1 \cap
  A_2\cap A_3) \cup ( A_1 \cap
  B_3)\cup B_2$, etc. Integrating this reasoning produces
  \begin{equation*}
    A_1(\phi)\subset
    \left(\bigcap_{i=1}^n 
      A_i(\phi)
    \right)\bigcup\left(\bigcup_{i=2}^n B_i(\phi)\right),
  \end{equation*} 
  so that
  \begin{equation*}
    A_1(\eps)\subset
    \left(\bigcap_{i=1}^n 
      A_i(\eps)
    \right)\bigcup\left(\bigcup_{i=2}^n B_i(\eps)\right)\;.
  \end{equation*} 
  Now define
  \begin{multline*}
    C_i(R)\eqdef\left\{\sup_{Q\in\Ha_{i-1}}\max_j\max_{k_1,k_2,\ldots,k_j}\sup_{\|\phi\|\leq
        1} \Norm{ \left\langle \bl
          F,Q,g_{k_1},\ldots,g_{k_j}\br(X(t)),K_{s,T}\phi\right\rangle}>R\right\}
    \\=\left\{\max_{Q\in\Extr\Ha_{i-1}}\max_j\max_{k_1,k_2,\ldots,k_j}\sup_{\|\phi\|\leq
        1} \Norm{ \left\langle \bl
          F,Q,g_{k_1},\ldots,g_{k_j}\br(X(t)),K_{s,T}\phi\right\rangle}>R\right\},
  \end{multline*}
  where $\Extr$ denotes the set of extreme points of a set. The second
  equality is implied by the fact that a linear function on a convex
  closed set attains its maximum at an extreme point of the set. Note
  also that $\Extr\Ha_i$ is finite for all $i$ (this can be proved by
  induction in $i$).
  
  Since Lemma~\ref{lm:induction_step_soft} implies $B_i(\eps) \subset
  H\left(\eps^{\kappa(i-1)}\right)\cup
  C_i\left(\eps^{-\kappa(i)}\right)$, we have
  \begin{equation*}
    A_1(\eps)\subset
    \bigcap_{i=1}^n 
    A_i(\eps)\cup\bigcup_{i=2}^n\left[ H\left(\eps^{\kappa(i-1)}\right)\cup
      C_i\left(\eps^{-\kappa(i)}\right)\right].
  \end{equation*} 
  Now setting 
  \begin{multline*}
    \widehat
    C_i(R)\eqdef\left\{\max_{Q\in\Extr\Ha_{i-1}}\max_j\max_{k_1,k_2,\ldots,k_j}
      \Norm{ \bl
        F,Q,g_{k_1},\ldots,g_{k_j}\br(X(t))}_\HH>R\right\}\;,
  \end{multline*}
  the second inequality in Lemma \ref{l:boundLipSupIP} implies that
  \begin{equation}\label{eq:splitCD}
    C_i(R) \subset \widehat C_i(\sqrt{R}/2) \cup D_*( \sqrt{R}/2)
  \end{equation}
(recall that $g_* \geq 1$).
  Defining
  \begin{xalignat*}{3}
    H_*(\eps)&=\bigcup_{i=2}^n
    H\left(\eps^{\kappa(i-1)}\right)&&\text{and} & C_*(\eps)&=\bigcup_{i=2}^n
    \widehat C_i\left(\eps^{-\kappa(i)/2}/2\right),
  \end{xalignat*}
  we have that 
 \begin{equation}\label{eq:partialDecomp}
   A_1(\eps)\subset\Bigl(\;
   \bigcap_{i=1}^n 
   A_i(\eps)\Bigr)\cup C_*(\eps) \cup H_*(\eps) \cup
   D_*(\eps^{-\kappa(n)/2}/2) . 
  \end{equation} 
  Observe that
  \begin{multline*}
    \bigcap_{i=1}^n A_i(\eps)\subset\left\{\inf_{\phi\in U}\sup_{Q\in \Ha_n}
      \supT\ip{Q(X(s))}{\phi} <\eps^{\kappa(n)}\right\} \\ \subset
    \left\{\inf_{\phi\in U}\sup_{Q\in
        \Ha_n}|\ip{Q(X(T))}{\phi}|<\eps^{\kappa(n)}\right\}\eqdef A_*(\eps).
  \end{multline*}
  Hence from \eqref{eq:partialDecomp} and the fact that $D_*(\eps^{-1/3})
  \subset D_*(\eps^{-\kappa(n)/2}/2)$, we have that
  \begin{multline}\label{eq:partialProbBound}
    \PP\left\{\inf_{\phi\in U}\ip{\MM(u(T))\phi}{\phi}<\eps\right\}
    \\\leq \PP\Big( D_*(\eps^{-\kappa(n)/2}/2) \cup A_*(\eps) \cup H_*(\eps)
    \cup C_*(\eps) \Big). 
  \end{multline}
  
  We now show that the probability of each of these terms is
  $o(\eps^q$) for any $q \geq 1$. Applying the Markov inequality and
  condition~\eqref{eq:approximation} yields
  \begin{equation*}
    \PP\big(A_*(\eps)\big) \leq  \Lambda_{q/\kappa(n)}^* \eps^q,
  \end{equation*}
  for all $q\geq 1$ and $\eps >0$. For all $\eps$ sufficiently
  small, the right hand side is less than $\eps^{q/2}$. 
  
  Lemma~\ref{lm:induction_step_soft} and the finiteness of $\Extr\Ha_i$
  imply that there are universal constants $K_1$, $K_2$ and $\gamma$,
  depending only on the number $n$ and the number of Brownian motions
  $d$, so that
  \begin{equation*}
    \PP\big( H_*(\eps)  \big) \leq K_1e^{-K_2\eps^{\gamma}}\;.
  \end{equation*}
  Turning to $D_*$, Lemma \ref{l:lipBoundsXK} implies that 
  \begin{align*}
    \PP\bigl(D_*(\eps^{-\kappa(n)/2}/2)\bigr)
     &\le \PP\left\{\Norm{K_{s,t}}> \eps^{-\kappa(n)/2}/2\right\}\\
     &\leq \left(C_{2q/\kappa(n)} \sqrt{K^*_{4q/\kappa(n)}\left(1+u^*_{4q(m-1)/\kappa(n)}\right)}\ 
      + K^*_{4q/\kappa(n)}\right)\eps^q\ . 
  \end{align*}
    Lastly from Lemma \ref{l:lipPoly}, Corollary \ref{l:unifLip}, and
  Assumption \ref{a:solReg}, we see that for any $q \geq 1$ there
  exists $p_1(n)\geq 1$ and a constant $C(n)$ so that
  \begin{equation*}
    \PP(C_*(\eps)) \leq \eps^q C(1+ u^*_{p_1})g_{*}^{p_1},
  \end{equation*}
  for any $\eps >0$. 
    
    Combining these bounds on the probability of the four sets with
  \eqref{eq:partialProbBound} completes the proof of the lemma.
\end{proof}

\begin{proof}[Proof of Lemma \ref{lm:induction_step_soft}]
  The proof begins the same way as that of Lemma \ref{lm:induction_step}.
  Upon reaching \eqref{eq:eqZeroLine}, we invoke Theorem~ \ref{th:chance_of_being_small_embedding_integral} 
  rather than Theorem~\ref{th:quadratic_variation_of_polynomial}.
\end{proof}

\section{Refinements and Generalizations}
\label{sec:refinments}
We now turn to a number of extensions and generalizations of the
preceeding results. In the first part of the section, we make more
explicit the dependence of the estimates on the initial data.
Understanding the dependence of the estimates on the initial data is
critical to proving results such
as unique ergodicity (see \cite{b:HairerMattingly04}). In the second half of
this section, we isolate the main arguments of this paper so that they
might be better applied to PDEs which do not fit into the precise setting of
this text.
\subsection{Dependence on the initial data}
\label{sec:depIntData}

\begin{theorem}\label{thm:refinedNonDegMal}
  In the setting of Section \ref{sc:framework}, assume that
  Assumptions \ref{a:solReg} and \ref{a:K2a} hold. Additionally,
  assume that there exists a function $\Psi:\HH \to
  [1,\infty)$ such that for any $p \geq
  1$ there exist constants $\tilde u_p^*(\TT,T)$ and $\tilde K_p^*(\TT,T)$
  so that
  \begin{align*}
     u_p^*(\TT,T,u_0) &\leq \tilde u_p^*(\TT,T) \Psi(u_0),\\
     K_p^*(\TT,T,u_0) &\leq \tilde K_p^*(\TT,T) \Psi(u_0),
  \end{align*}
for all $u_0 \in \HH_0$.

Consider the setting of Theorem \ref{lm:moments_matrix_inverse}. If either 
\begin{enumerate}
\item $S$ is a finite-dimensional subset of $\Hbb_n$ for some
  $n< \infty$, or
\item $S$ is a finite-dimensional subset of $\VV$ so that for some $n<
  \infty$ and for any $p\geq 1$, the condition given in equation
  \eqref{eq:approximation} holds. Furthermore, for any $p \geq 1$, there
  exists a positive constant $\tilde\Lambda_p^*(\TT,T)$ so that
  \begin{align*}
      \Lambda_p^*(\TT,T,u_0) &\leq \tilde \Lambda_p^*(\TT,T) \Psi(u_0),
  \end{align*}
  where $\Lambda_p^*$ was also defined in \eqref{eq:approximation},
\end{enumerate}
then for  any $p \geq 1$, there are positive constants
$C$, $\eps_0$, $q$, and $\delta$  such that   
  \begin{equation}\label{eq:evProp}
    \PP\left\{\inf_{\phi\in
        U_\delta}\ip{\MM(u(T))\phi}{\phi}<\eps\right\}\leq C \Psi^q(u_0)\eps^p,   
  \end{equation}
  for all $\eps \in (0,\eps_0]$ and $u_0 \in \HH_0$. Here, as
  before, $U_\delta=\{\phi\in \HH: \|\phi\|\le 1,\ \|\Pi\phi\|\geq \delta \},$ 
  and $\Pi$ is
  the projection onto $S$. In the first case
  $C$ depends on $p,\TT,T,S,u_p^*$, and $K_p^*$, and in the second it also depends
  on $\Lambda_p^*$. In both cases  $\eps_0$ depends only on $S$,
  and $q$ depends only on $p$ and $S$.

\end{theorem}

\begin{proof}
  Looking back at the proof of Theorem
  \ref{lm:moments_matrix_inverse}, we need to obtain a bound of the
  quoted type on the right hand side of \eqref{eq:partialProbBound}.
  In light of the calculations in the proof bounding the size of the
  various sets, the probabilities of $D_*$, $\widehat C_*$, and $H_*$ are
  all bounded as desired because of the assumptions of Theorem
  \ref{thm:refinedNonDegMal}. The only set left uncontrolled is
  $A_*$.

However, all the vector fields in $\Hbb_n$ are constant, and hence there
is an $\eps_0$ sufficiently small and depending only on the
structure and size of $\Hbb_n$ and the $S$ chosen so that, if $\eps
\in (0,\eps_0]$, then $A(\eps)$ is empty.
\end{proof}

\subsection{Generalizations}

We now state a few ``meta'' theorems. The assumptions require extra
work to verify but they isolate the main parts of the argument and
allow the ideas to be applied to a wider range of PDEs which do not fit
exactly into the previous settings. We relax our assumptions on $N$,
assuming only that it is a polynomial from  $\Dom(L)$ into $\HH$. We
assume that with probability one
\begin{equation*}
  u\in C\big([0,T],\HH\big)\cap  L^\infty_{loc}\big((0,T],\Dom(L)\big).
\end{equation*}
Lastly we fix a Banach space  $(\HH_1,\lvert\;\cdot\;\lvert_{\HH_1})$, with $\HH_1 \subset \HH$, and assume that 
for each $g_k$ and $\phi \in \HH_1$
\begin{equation*}
  \ip{g_k}{K_{t,T}\phi} \in C\big([\Tt,T],\RR\big)
\end{equation*}
with probability one as a function of $t$. We now define a new set of
admissible vector fields.
\begin{definition}
  $\widetilde \Ab$ is  the set of all
  polynomial vector fields $Q:\VV\to\VV'$ such that with probability one the
  following conditions hold:
  \begin{enumerate}
  \item $Q(X(t))\in L^2\big([\Tt,T],\VV\big)$,
  \item $\frac{d}{dt}Q(X(t))\in L^2\big([\Tt,T],\VV'\big)$,
  \item For all $0\leq i\leq m$, $k_j \in\{ 1,\cdots,m\}$ and $\phi \in \HH_1$, 
    \begin{equation}
      \label{eq:XinDom}
      \ip{[F,Q,g_{k_1},\dots,g_{k_i}](X(t))}{K_{t,T}\phi} 
    \end{equation}
    is well defined and in $C\big([\Tt,T],\RR\big)$ as a function of $t$.
 \end{enumerate}
\end{definition}
 
Next, define $\widetilde \Hb_n$ exactly as in \eqref{a:HaDef}, replacing
$\Ab$ by $\widetilde\Ab$.

\begin{theorem}\label{thm:genABS}  Assume that Assumptions \ref{a:J1},
  \ref{a:K1}, and \ref{a:D2}  hold.
  Let $S$ be a finite-dimensional linear subspace which is a subset of
  $\widetilde \Hb_n(X(T))\cap \HH_1$ with probability one. Then the
  distribution of the 
  projection of $X(T)$ onto $S$ is absolutely continuous with respect
  to Lebesgue measure on $S$. 
\end{theorem}

Turning to smoothness, define $\mathcal{L}^\infty$ to be the space of all processes $f:[0,T] \to
\RR$ such that $\E \Norm{f}^p < \infty$, for all $p \geq 1$. Let
$\widetilde \Aa$ be defined by
\begin{multline*}
  \widetilde \Aa\eqdef \Big\{ Q \in \widetilde \Ab :\sup_{\substack{\phi \in
    \HH_1\\ \|\phi\|\le 1 }}\Norm{\ip{[F,Q,g_{k_1},\dots,g_{k_i}](X(t))}{K_{t,T}\phi}}
  \in \mathcal{L}^\infty, \\ \text{ for all } 0\leq i\leq m \text{ and
  } k_j \in\{ 1,\cdots,m\} \Big\}.
\end{multline*}
Lastly, define $\widetilde \Ha_n$ as in \eqref{eq:HaDef}, but with
$\Aa$ replaced by $\widetilde \Aa$.

\begin{theorem}\label{th:smoothnessC}
  Assume that Assumptions \ref{a:J1} and \ref{a:K1} 
  hold.  Let $S$ be a deterministic finite-dimensional subspace of
  $\HH_1$ such that for some $n$ and $\delta>0$,
\begin{equation}
  \label{eq:approximationGen}
\widetilde \Lambda_p^*(u_0,T) \eqdef \E\left[\inf_{\substack{\phi\in \tilde U_\delta}
        }\sup_{Q\in \widetilde\Ha_n(X(T))}\Bigl|\langle
    \phi,Q(X(T))\rangle\Bigr|\right]^{-p}<\infty, 
\end{equation}
for all $p\geq1$. Here  $\tilde U_\delta=\{\phi\in \HH_1: |\phi|_{\HH_1}\le1,\ |\Pi\phi|_{\HH_1}\geq \delta \}$.
If $\Pi_S u(T) \in \D^\infty_S$, then the density of $\Pi_S u(T)$ with
respect to Lebesgue measure (whose existence is guaranteed by
Theorem~\ref{thm:genABS}) is a $C^\infty$-function on $S$.
\end{theorem}

In the spirit of section \ref{sec:depIntData}, we now give a ``meta''
theorem which isolates the dependence on the initial data.

\begin{theorem} As above, assume that Assumptions \ref{a:J1} and
  \ref{a:K1} hold.  Let $S$ be a deterministic finite-dimensional
  subspace of $\HH_1$ such that, for some $n$ and $\delta$, the bound in
  \eqref{eq:approximationGen} holds.

  Let $\Psi:\HH_0 \rightarrow (0, \infty)$ be a function such that,  for
    all $p \geq 1$, there exists a $C_p$ such that:
  \begin{enumerate}
  \item For any $Q \in \widetilde \Ha_n$,
    \begin{equation*}
      \E \sup_{\substack{\phi \in \HH_1\\ |\phi|_{\HH_1}\le 1}}
      \Norm{\ip{[F,Q,g_{k_1},\dots,g_{k_i}](X(t))}{K_{t,T}\phi}}^p  \leq
     C_p \Psi(u_0),
   \end{equation*}
   for all $u_0 \in \HH_0$, $0\leq i\leq m$, and $k_j \in\{
   1,\cdots,m\}$;
 \item $\widetilde \Lambda_p^*(u_0,T) \leq C_p \Psi(u_0)$.
 \end{enumerate}
 Then the conclusion given in \eqref{eq:evProp} holds with $U$
 replaced by the $\tilde U_\delta$ defined in Theorem \ref{th:smoothnessC}
 and for constants with the same dependencies as in
 Theorem \ref{thm:refinedNonDegMal}.
\end{theorem} 

\section {Non-adapted polynomials of Wiener processes}
\label{sc:non-adapted_stuff}
This section contains the technical estimates which are the heart of
the paper. They are the key steps in the proofs in
Section~\ref{sc:genResults} which ensure that the randomness moves,
with probability one, to all of the degrees of freedom connected to
the noise directions through the nonlinearity.  The results in
section~\ref{sc:non-adapted_qualitative} are more qualitative and are
the basis of the proof of existence of absolutely continuous
densities.  Section~\ref{sc:more_quantative_estimates} contains the
more quantitative estimates needed to prove the smoothness of the
density and give estimates on the eigenvalues of the Malliavin matrix.
That being said, the basic ideas of the two sections are the same. We
show that coefficients of a finite Wiener polynomial (see below for more
details) are small with high probability if the entire polynomial is
small, even if the coefficients are not adapted to the Wiener
processes.

The core idea, used in our context, dates back at least to the
pioneering work of Malliavin, Bismut, Stroock and others on the
probabilistic proof of the existence of smooth densities for hypoelliptic
diffusions in finite dimensions. The techniques developed there (see
\cite{b:KusuokaStroock84,b:Norris86}) used martingale estimates to
relate the size of a process to its quadratic variation.  Here we 
cannot make use of such martingale estimates directly since we have
non-adapted stochastic processes. The non-adaptedness arose in a
natural way because we only have a semiflow and cannot return all
estimates to the tangent space at the origin and work with the reduced
Malliavin covariance matrix which is adapted.  As is often done, we
replace an adaptedness assumption with an assumption on the regularity
in time of the processes. This section is a generalization of the
results in \cite{MattinglyPardoux05} which proved similar results
for quadratic polynomials of Wiener processes. The proofs here extend
these results to polynomials of any order while also simplifying the
proofs.

\subsection{Qualitative results}
\label{sc:non-adapted_qualitative}

Consider a probability space $(\Omega,\Fc,\PP)$. For a stochastic
process $X$ defined on $[0,t]$, we define $\Delta_{s_1s_2}(X) =
X(s_2)-X(s_1)$ .  For two stochastic processes $X_1,X_2$ defined on
the same time interval $I=[\Ta,\Tb]$, we denote
\begin{equation*}
  \langle X_1,X_2\rangle_I\eqdef\lim_{N\to\infty} \sum_{j=1}^N\Delta_{t_{j-1}t_j}(X_1)
  \Delta_{t_{j-1}t_j}(X_2)\quad\mbox{ in probability},
\end{equation*}
if this limit exists, where $\Ta=t_0^N<\ldots<t_N^N=\Tb$ for each $N$
and $\sup\{t_j^N-t_{j-1}^N\}\to 0$ as $N\to\infty$. We shall also
write $\langle X\rangle_I=\langle X,X\rangle_I$ and $\langle
X\rangle_t=\langle X\rangle_{[0,t]}$.

We begin by considering the basic cross quadratic variation between
two monomial terms. We emphasize that the processes $A(s)$ and $B(s)$ in the
following lemma need not be adapted to the filtration generated by the
Wiener processes. 
\begin{theorem}
  \label{th:quadratic_variation}
  Let $W_1(s),W_2(s),\ldots, W_d(s)$ be a collection of mutually
  independent standard one-dimensional Brownian motions on a time
  interval $I$ and let $A(s), B(s)$ be two continuous and bounded
  variation stochastic processes defined on $I$.  Then
  \begin{multline*}
    \langle AW_{i_1}\ldots W_{i_n},BW_{k_1}\ldots W_{k_m}\rangle_I \\
    =\int_I\left[A(s)B(s)\sum_{p=1}^n\sum_{q=1}^m \delta_{i_pk_q}
      \frac{W_{i_1}(s)\ldots W_{i_n}(s)W_{k_1}(s)\ldots
        W_{k_m}(s)}{W_{i_p}(s)W_{k_q}(s)}\right]ds.
  \end{multline*}
\end{theorem}

\bpf In the proof we write $W_i(j)$ instead of $W_i(t_j^N)$, $A(j)$
instead of $A(t_j)$,
$\Delta_{i,j}$ instead of $\Delta_{t_i,t_j}$, and $j^-$ instead of
$j-1$. We begin by observing that
\begin{equation}
  \label{eq:sum in terms of Q}
  \sum_{j=1}^N\Delta_{j^-,j} (AW_{i_1}\ldots W_{i_n})
  \Delta_{j^-,j}(BW_{k_1}\ldots W_{k_n})=\sum_{j=1}^N Q_j^{(A)} Q_j^{(B)},
\end{equation}
where
\begin{multline*}
  Q_j^{(A)}=\Delta_{j^-,j}(A) W_{i_1}(j^-)\ldots W_{i_n}(j^-) \\+A(j)
  \Delta_{j^-,j} (W_{i_1}) W_{i_2}(j^-)\ldots W_{i_n}({j^-})
  \\+A(j)W_{i_1}({j})\Delta_{j^-,j}(W_{i_2})W_{i_3}({j^-})\ldots
  W_{i_n}({j^-}) +\ldots \\+A(j)W_{i_1}({j})\ldots
  W_{i_{n-1}}(j)\Delta_{{j^-},j} (W_{i_n}),
\end{multline*}
and
\begin{multline*}
  Q_j^{(B)}=\Delta_{{j^-},j}(B) W_{k_1}({j^-})\ldots W_{k_m}({j^-})\\
  +B(j) \Delta_{{j^-},j} (W_{k_1}) W_{k_2}({j^-})\ldots W_{k_m}({j^-})\\
  +B(j)W_{k_1}({j})\Delta_{{j^-},j}(W_{k_2})W_{k_3}({j^-})\ldots
  W_{k_m}({j^-}) +\ldots \\+B(j)W_{k_1}({j})\ldots
  W_{k_{m-1}}(j)\Delta_{{j^-},j} (W_{k_m}).
\end{multline*}

Therefore, the sum in \eqref{eq:sum in terms of Q} contains the following terms:

\begin{equation*}
  \sum_{j=1}^N \Delta_{{j^-},j}(A) W_{i_1}({j^-})\ldots
  W_{i_n}({j^-})\Delta_{{j^-},j}(B) W_{k_1}({j^-})\ldots
  W_{k_m}({j^-}), 
\end{equation*}
\begin{equation*}
  \sum_{j=1}^N \Delta_{{j^-},j}(A) W_{i_1}({j^-})\ldots
  W_{i_n}({j^-})B(j) W_{k_1}({j^-}) 
  \ldots \Delta_{{j^-},j}(W_{k_q})\ldots W_{k_m}({j^-}), 
\end{equation*}
\begin{equation*}
  \sum_{j=1}^N A(j) W_{i_1}({j^-})\ldots \Delta_{{j^-},j}(W_{k_p})
  \ldots W_{i_n}({j^-})\Delta_{{j^-},j}(B) W_{k_1}({j^-})\ldots
  W_{k_m}({j^-}), 
\end{equation*}
\begin{equation*}
  \sum_{j=1}^N A(j) W_{i_1}({j^-})\ldots \Delta_{{j^-},j}(W_{k_p})
  \ldots W_{i_n}({j^-}) B W_{k_1}({j^-})\ldots
  \Delta_{{j^-},j}(W_{k_q})\ldots  W_{k_m}({j^-}).
\end{equation*}

The first three sums above converge to zero as $t_j-t_{j-1}\to 0$,
since $A$ and $B$ are of bounded variation and continuous and all
$W_i$ are continuous. Lemmas~4.2 and~4.3 from \cite{MattinglyPardoux05} imply that the
fourth sum above converges to 
\begin{equation*}
  \int_I\left[A(s)B(s) \delta_{i_pk_q} \frac{W_{i_1}(s) \ldots
      W_{i_n}(s)W_{k_1}(s) \ldots W_{k_m}(s)}{W_{i_p}(s)W_{k_q}(s)}\right]ds, 
\end{equation*}
and the theorem is proved.\epf

\begin{corollary}\label{c:Zqual}
  Let $\Ac$ be collection of stochastic processes on $\Omega$
  such that there is a set $\Omega'\in\Fc$, with $\PP(\Omega')=1$,
  so that for each $\omega\in\Omega'$ all of the process in $\Ac$ are
  of bounded variation and continuous.

  Then there is a set $\Omega''\subset\Omega'$, with
  $\PP(\Omega'')=1$, and a sequence of partitions
  $t^{(N)}=\{\Ta=\Tt^N<\ldots<t_N^N=\Tb\}$, with
  $\sup\{t_j^N-t_{j-1}^N\}\to 0$, as $N\to\infty$, such that for any
  process $Z(t)$ of the form
 \begin{equation*}
    Z=A^{(0)}+\sum_{i_1}A_{i_1}^{(1)}W_{i_1}+
    \sum_{i_1,i_2}A_{i_1,i_2}^{(2)}W_{i_1}W_{i_2} + \ldots
    +\sum_{i_1,\ldots,i_n}A^{(n)}_{i_1,\ldots,i_n}W_{i_1} \ldots W_{i_n},
  \end{equation*}
  with  $A_{i_1,\cdot,i_k}^{(k)} \in \Ac$, one has that the limit
\begin{equation*}
  \lim_{N\to\infty} \sum_{j=1}^N\Delta_{t_{j-1}t_j}^2(Z)     
\end{equation*}
exists on $\Omega''$ and equals $\langle Z\rangle_I$.
\end{corollary}
\begin{proof} We notice that the proof of Theorem~\ref{th:quadratic_variation} 
implies that there is a full measure set $\tilde\Omega$  that 
is defined in terms of the Wiener processes involved, with the following property:
if for $\omega\in\tilde\Omega$ the realization of a process $A$ possesses the mentioned regularity properties,
then the desired convergence holds. The proof is
completed by setting $\Omega''=\Omega'\cap\tilde\Omega$.
\end{proof}

We now use the previous results to prove that in the setting of the
previous corollary, if $Z$ is identically zero, then the coefficients
$A_{i_1,\cdot,i_k}^{(k)}$ must be identically zero.
\begin{theorem}\label{th:quadratic_variation_of_polynomial}
  Let $\Ac$ and $Z$ be as in the above corollary, and let $\Omega''$
  be the set given in the conclusion of the same corollary. Additionally, assume
  that, for each $\alpha$ and $i_1,\ldots,i_\alpha\in 1,\ldots,d$, the
  coefficients $A^{(\alpha)}_{i_1,\ldots,i_\alpha}$ are symmetric
  (i.e. invariant under substitutions on indices
  $i_1,\ldots,i_\alpha$). 

  If $Z(s)=0$ for all $s\in[0,T]$ with probability one, then all
  the processes $A^{(\alpha)}_{i_1,\ldots,i_\alpha}$ are identically
  zero on $[0,T]$ with probability one.
\end{theorem}
\begin{proof}
 We proceed by induction. For $n=0$ the statement of the theorem
is obvious. Now suppose $n\geq1$. Then
\begin{multline*}
  \langle
  Z\rangle_T=\sum_{\alpha,\beta=1}^n\left\langle\sum_{i_1,\ldots,i_\alpha}A^{(\alpha)}_{i_1,\ldots,i_\alpha}W_{i_1}\ldots  
    W_{i_\alpha},
    \sum_{k_1,\ldots,k_\beta}A^{(\beta)}_{k_1,\ldots,k_\beta}W_{k_1}\ldots
    W_{k_\beta}\right\rangle_T
  \\=\int_0^T\sum_{\alpha,\beta=1}^n\sum_{i_1,\ldots,i_\alpha}\sum_{k_1,\ldots,k_\beta}A^{(\alpha)}_{i_1,\ldots,i_\alpha}A^{(\beta)}_{k_1,\ldots,k_\beta}
  \sum_{p=1}^\alpha\sum_{q=1}^\beta \delta_{i_pk_q}\frac{W_{i_1}\ldots
    W_{i_\alpha}}{W_{i_p}}\frac{W_{k_1}\ldots W_{k_\beta}}{W_{k_q}}ds
  \\=\int_0^T\sum_{r=1}^d\sum_{\alpha,\beta=1}^n\sum_{i_1,\ldots,i_\alpha}\sum_{k_1,\ldots,k_\beta}
  \sum_{p=1}^\alpha\sum_{q=1}^\beta\delta_{ri_p}\delta_{ri_q}
  \frac{A^{(\alpha)}_{i_1,\ldots,i_\alpha}W_{i_1}\ldots
    W_{i_\alpha}}{W_{i_p}} \frac{A^{(\beta)}_{k_1,\ldots,k_\beta}
    W_{k_1}\ldots W_{k_\beta}}{W_{k_q}} ds
  \\=\sum_{r=1}^d\int_0^T\left(\sum_{\alpha=1}^n\sum_{i_1,\ldots,i_\alpha}
    \sum_{p=1}^\alpha\delta_{ri_p}
    \frac{A^{(\alpha)}_{i_1,\ldots,i_\alpha}W_{i_1} \ldots
      W_{i_\alpha}}{W_{i_p}}\right)^2 ds.
\end{multline*}
Since we assumed that $Z(s)=0$ for $s\in[0,T]$ and the integrand is
continuous, we conclude that
\begin{equation*} 
  Z_r(s)\eqdef\sum_{\alpha=1}^n\sum_{i_1,\ldots,i_\alpha} \sum_{p=1}^\alpha\delta_{ri_p} \frac{A^{(\alpha)}_{i_1,\ldots,i_\alpha}(s)W_{i_1}(s)\ldots
    W_{i_\alpha}(s)}{W_{i_p}(s)}=0
\end{equation*}
for each $r=1,\ldots,d$ and all $s\in[0,T]$. Notice now that due to
the symmetry of coefficients $A$ the process $Z_r(s)$ satisfies the
assumptions of the theorem with $n$ reduced by one.  That
$A^{(\alpha)}_{i_1,\ldots,i_\alpha}(s)=0$ a.s. for $\alpha\geq1$,
follows from the fact that all coefficients of $Z_{i_1}(s)$ are equal
to zero a.s. by the induction hypothesis. Since $Z\equiv 0$ and
$A^{(\alpha)}\equiv0$ for positive $\alpha$, we conclude that
$A^{(0)}\equiv0$ as well. The theorem is proved.
\end{proof}

\subsection{More Quantitative Estimates}\label{sc:more_quantative_estimates}

Now our aim is to prove a quantitative version of the last
theorem. Again we consider a process $Z(t)$ of the same form as in
Corollary \ref{c:Zqual}. To do so we introduce a family of Wiener
polynomials with constant coefficients which will be used to
approximate $Z$. Namely, for any nonnegative integer $n$ and collection
of coefficients $\lambda$ with
\begin{equation*}
  \lambda=\left\{\lambda^{(\alpha)}_{i_1,\ldots, i_\alpha} \in
    \RR,\alpha=0,\ldots,n, i_1,\ldots, i_\alpha=1,\ldots,d\right\}, 
\end{equation*}
we define
\begin{equation*}
  Z_\lambda=\lambda^{(0)}+\sum_{i_1}\lambda_{i_1}^{(1)}W_{i_1}+
    \sum_{i_1,i_2}\lambda_{i_1,i_2}^{(2)}W_{i_1}W_{i_2}+\ldots
    +\sum_{i_1,\ldots,i_n}\lambda^{(n)}_{i_1,\ldots,i_n}W_{i_1}\ldots W_{i_n}.
\end{equation*}
We now introduce a collection of typical coefficients, a set of typical
Wiener processes, and a collection of atypical $Z_\lambda$, which are
too small in light of their coefficients not being uniformly small.
This last set captures the event which we wish to describe, but for
the $Z_\lambda$ rather than the $Z$. We begin with the coefficients
$\lambda$, which we do not want to be uniformly too small.

For a real number $\eps>0$ and a nonnegative integer $n$ define
$\Lambda(\eps,n)$ to be the set of coefficients
$\lambda=\left\{\lambda^{(\alpha)}_{i_1,\ldots, i_\alpha},\alpha=0,\ldots,n,
  i_1,\ldots, i_\alpha=1,\ldots,d\right\}$ such that
\begin{equation*}
  \max \left\{|\lambda^{(\alpha)}_{i_1,\ldots,
      i_{\alpha}}|:\alpha=0,\ldots,n, i_1,\ldots, i_\alpha=
    1,\ldots,d\right\}\geq\eps. 
\end{equation*}
We now define a set of atypical $Z_\lambda$ with $\lambda \in
\Lambda(\eps,n)$.  Take $\hat \eps>0$ and divide the segment $[0,T]$
into $m=[T \hat \eps^{-\frac32 8^{n+1}}]+1$ segments
$I_1=[0,t_1],I_2=[t_1,t_2],\ldots,I_{m}=[t_{m-1},t_m]$, each one of
length less than $\hat \eps^{\frac32 8^{n+1}}$ and greater than
$\frac12\hat \eps^{\frac32 8^{n+1}}$.

Let
\begin{equation*}
  D^*(\hat\eps,I,\Lambda(\eps,n))=\left\{\inf_{\lambda\in\Lambda(\eps,n)}\sup_{t\in
    I}|Z_\lambda(t)|<\hat\eps\right\}
\end{equation*}
and define
\begin{equation*}
  F(\hat\eps,\eps)=\bigcup_{k=1}^m D^*\left(\hat\eps,I_k,\Lambda(\eps,n)\right).  
\end{equation*}

To define the set of typical Wiener trajectories, recall that for any
function $f:[0,T]\to\R$ we define its $\rho$-H\"older constant by
\begin{equation*}
  \Hol_{\rho}(f)\eqdef\sup_{\substack{0\leq s<r \leq T}}
\frac{|f(s)-f(r)|}{|s-r|^\rho},
\end{equation*}
and  
\begin{equation*}
  \Norm{f}_\rho\eqdef\max\{\norm{f}_{L^\infty},\Hol_{\rho}(f)\}.
\end{equation*}
With this definition, we introduce the set of Wiener processes
\begin{equation*}
  B(R)=\left\{\Norm{W_{i_1}\ldots W_{i_\alpha}}_{\frac14}<R, \alpha=1,\ldots,n, 
i_1,\ldots,i_\alpha=1,\ldots,d\right\}.
\end{equation*}
\begin{remark}
  Notice that the sets $B$ and $F$ are universal in that they do not
  depend on the processes $A$ in any way other than through the number $n$. 
\end{remark}

We now are ready to state the quantitative version of Corollary
\ref{c:Zqual}. We want to conclude that if $Z$ is small it is unlikely
that the $A$ processes are not small. The sets $D$ and $E$ below embody the
first event and the complement of the second event, respectively:
 \begin{align*}
  D(\eps)&=\{\norm{Z}_{L^\infty}<\eps\},\\
  E(\eps)&=\left\{\max_{\alpha=1,\ldots,n}
    \max_{i_1,\ldots,i_\alpha=1,\ldots,d}
    \norm{A^{(\alpha)}_{i_1,\ldots, i_\alpha}}_{L^\infty}<
    \eps\right\}.
\end{align*}
To state the result we need to define a localization set which ensures
that we can well approximate $Z$ by a $Z_\lambda$ process with
$\lambda \in \Lambda(\eps,n)$. Defining
\begin{equation*}
  C(R)=\left\{\Lip(A^{(\alpha)}_{i_1,\ldots, i_\alpha})<R, \alpha=1\ldots,n, 
i_1,\ldots,i_\alpha=1,\ldots,d\right\},
\end{equation*}
we have the desired results. 
\begin{theorem}\label{th:chance_of_being_small_embedding} 
For each $n$ there is $\eps_0(n)$ 
depending only on $n,d$ and $T$ such that
\begin{equation}
D(\eps^{8^{n+2}})\cap E^c(\eps)\cap  C(\eps^{-1}) \subset  
B^c(\eps^{-1/5})\cup F(\eps^{8^{n+1}\cdot 5/4 },\eps^{2+\frac{1}{n+1}}),
\label{eq:inclusion_in_unlikely_set}
\end{equation}
for all $\eps<\eps_0.$ 
\end{theorem}

\begin{theorem}\label{th:chance_of_being_small_estimate} 
  For each $n$ there are positive numbers $\eps_1(n), q_1(n),K_1(n),K_2(n)$ 
  depending only on $n,d$ and $T$ such that
  \begin{equation*}
    \PP\left(B^c(\eps^{-1/5})\cup 
      F(\eps^{8^{n+1}\cdot 5/4},\eps^{2+\frac{1}{n+1}})\right)
    <K_1\exp\{-K_2\eps^{-q_1}\},
  \end{equation*}
  if $\eps<\eps_1.$
\end{theorem}

\begin{remark}
  Theorem~\ref{th:chance_of_being_small_estimate} provides an estimate
  of the set appearing in the statement of
  Theorem~\ref{th:chance_of_being_small_embedding}.  Thus, these two
  theorems say that if $Z$ is small ( the event $D(\eps^{8^{n+2}})$), then
  with high probability the coefficients $A$ defining $Z$ are small as
  well (the event $E(\eps)$) on the localization set $C(\eps^{-1})$.
  Since the  $A$ are not necessarily adapted, one aim of
  Theorem~\ref{th:chance_of_being_small_embedding} is to reduce the
  problem to the traditional stochastic It\^o calculus.
  Notice also that the events in the r.h.s.
  of~\eqref{eq:inclusion_in_unlikely_set} are defined only in terms of
  the Wiener processes $W$.   
\end{remark}

We will in fact find not that $Z$ is uniformly small in time, but
rather that its integral in time is small.  However the following
results show how to reduce this case to the previously considered
setting.
 
Consider an arbitrary $\R$-valued random variable $g_0$ and define
\begin{xalignat*}{3}
  g(t)&=g_0+\int_0^tZ(s)ds,&
  \overline D(\eps) &= \left\{\|g\|_{L^\infty}<\eps\right\},&
\overline C(R) &= C(R)\cap E(R).
\end{xalignat*}

\begin{theorem}\label{th:chance_of_being_small_embedding_integral} 
  For each $n$ there is $\eps_0(n)$ depending only on $n,d$ and $T$
  such that
\begin{equation}
  \overline D(\eps)\cap E^c(\eps^{8^{-(n+3)}})\cap  \overline
  C(\eps^{-8^{-(n+3)}}) \subset   
  B^c(\eps^{-8^{-(n+3)}/5})\cup
  F(\eps^{\frac{5}{256}},\eps^{8^{-(n+3)}\left(2+\frac{1}{n+1}\right)}), 
\label{eq:inclusion_in_unlikely_set_integral}
\end{equation}
for all $\eps<\eps_0.$
\end{theorem}

The probability of the r.h.s. is estimated in the following theorem, which
is a direct consequence of
Theorem~\ref{th:chance_of_being_small_estimate}:
\begin{theorem}\label{th:chance_of_being_small_estimate_integral} 
  For each $n$ and numbers $\eps_1(n), q_1(n),K_1(n),K_2(n)$ 
  defined in Theorem~\ref{th:chance_of_being_small_estimate},
  \begin{equation*}
    \PP\Big( B^c(\eps^{-8^{-(n+3)}/5})\cup
    F(\eps^{\frac{5}{256}},\eps^{8^{-(n+3)}\left(2+\frac{1}{n+1}\right)}) \Big)
    <K_1\exp\{-K_2\eps^{-8^{-(n+3)}q_1}\}
  \end{equation*}
  if $\eps<\eps_1^{-8^{-(n+3)}}.$
\end{theorem}

Theorem~\ref{th:chance_of_being_small_embedding_integral} will follow
from Theorem~\ref{th:chance_of_being_small_embedding} and the next
lemma taken from \cite{MattinglyPardoux05}. We will give the
proof of Theorem~\ref{th:chance_of_being_small_embedding_integral}
before returning to the proof of
Theorem~\ref{th:chance_of_being_small_embedding} and
Theorem~\ref{th:chance_of_being_small_estimate}.

\begin{lemma}\label{l:MP2}\cite[Lemma 7.4]{MattinglyPardoux05} Let
  \begin{equation*}
    G(t)=G_0+\int_0^tH(s)ds,
  \end{equation*}
  where $G$ and $H$ are $\R$-valued functions and $G_0\in\R$. Suppose
  $\Hol_\alpha(H)\leq c\eps^{-\gamma}$ for some fixed $\alpha>\gamma>0$ and $\eps>0$.
  If $t\geq\eps^{\frac{1+\gamma}{1+\alpha}}$, then $\|G\|_\infty\leq\eps$ implies
  $\|H\|_\infty\leq (2+c)\eps^{\frac{\alpha-\gamma}{1+\alpha}}.$
\end{lemma}
\bpf[Proof of Theorem~\ref{th:chance_of_being_small_embedding_integral}] 

We begin by considering a generic term $A^{(k)}_{i_1,\ldots,i_k}W_{i_1}\ldots W_{i_k}$ from $Z$. On $
B(\eps^{-8^{-(n+3)}/5})\cap \overline C(\eps^{-8^{-(n+3)}})$, we have that 
\begin{align*}
  \Hol_{1/4}(A^{(k)}_{i_1,\ldots,i_k}W_{i_1}\ldots W_{i_k}) \leq& \Lip(A^{(k)}_{i_1,\ldots,i_k})\|W_{i_1}\ldots W_{i_k}\|_{L^\infty}\\ &\qquad \qquad +  \Hol_{1/4}(W_{i_1}\ldots W_{i_k})\|A^{(k)}_{i_1,\ldots,i_k}\|_{L^\infty}\\
\leq& 2\eps^{-2\cdot 8^{-(n+3)}}.
\end{align*}
Since there are no more than $d^n$ such terms for each degree between
$0$ and $n+1$, on $ B(\eps^{-8^{-(n+3)}/5})\cap \overline
C(\eps^{-8^{-(n+3)}})$ we have
\begin{equation*}
  \Hol_{1/4}(Z)
  <2(n+1)d^n\eps^{-2\cdot 8^{-(n+3)}}.
\end{equation*}
Then Lemma \ref{l:MP2} implies
\begin{equation}
  \label{eq:Z_in_terms_of_eps}
  \|Z\|_{\infty}\leq(2+2(n+1)d^n)
  \eps^{\frac{1/4-2\cdot 8^{-(n+3)}}{1+1/4}}.
 \end{equation}
Define 
\begin{equation*}
  \delta=\eps^{8^{-(n+3)}}.  
\end{equation*}
Then \eqref{eq:Z_in_terms_of_eps} implies that for small $\eps$ on
$\overline D(\eps)\cap B(\eps^{-8^{-(n+3)}})\cap \overline C(\eps^{-8^{-(n+3)}})$
\begin{equation*}
  \|Z\|_{\infty}\leq\delta^{ \frac{1/4-2\cdot 8^{-(n+3)}}{(1+1/4)8^{-(n+3)}} }
  \leq \delta^{8^{(n+2)}},  
\end{equation*}
\textit{i.e.}
\begin{equation}
  \label{eq:overlineD_inclusion}
  \overline D(\eps)\cap B(\eps^{-8^{-(n+3)}/5})\cap \overline C(\eps^{-8^{-(n+3)}})
  \subset D(\delta^{8^{(n+2)}}).
\end{equation}
Next,
\begin{align}
  \label{eq:main_inclusion_for_overlineD}
  \overline D(\eps)&\cap B(\eps^{-8^{-(n+3)}/5})\cap \overline
  C(\eps^{-8^{-(n+3)}}) \cap E^c(\delta)\notag \\&=\overline D(\eps)\cap
  B(\eps^{-8^{-(n+3)}/5})\cap \overline C(\eps^{-8^{-(n+3)}}) \cap
  E^c(\delta)\cap C(\delta^{-1})\notag \\&\subset
  D(\delta^{8^{(n+2)}})\cap E^c(\delta)\cap C(\delta^{-1})\cap
  B(\delta^{-1/5})\notag \\&\subset
  F(\delta^{8^{n+1}\cdot5/4},\delta^{2+\frac{1}{n+1}}),
\end{align}
where the identity is implied by $\overline
C(\eps^{-8^{-(n+3)}})\subset C(\delta^{-1})$, the first inclusion is a
consequence of~\eqref{eq:overlineD_inclusion} and
$B(\delta^{-1/5})=B(\eps^{-8^{-(n+3)}/5})$, and the second one from
Theorem~\ref{th:chance_of_being_small_embedding}.  Now
\eqref{eq:inclusion_in_unlikely_set_integral} is equivalent
to~\eqref{eq:main_inclusion_for_overlineD}, and the proof is complete.
\epf

We now return to the proofs of the central results of this section.
\begin{proof}[Proof of Theorem~\ref{th:chance_of_being_small_embedding}]  
Consider
\begin{equation*}
  G(\eps)=D(\eps^{8^{n+2}})\cap E^c(\eps)\cap  C(\eps^{-1}) \cap B(\eps^{-1/5}).  
\end{equation*}
To prove the theorem
it is sufficient to show
\begin{equation*}
  G(\eps)\subset F(\eps^{8^{n+1}\cdot5/4},\eps^{2+\frac{1}{n+1}}).  
\end{equation*}
We have
\begin{equation*}
  G(\eps)\subset\bigcup_{k=1}^m G_k(\eps),
\end{equation*}
where
\begin{equation*}
  G_k(\eps)=D(\eps^{8^{n+2}},I_k)\cap E^c(\eps,I_k)\cap  C(\eps^{-1}) \cap B(\eps^{-1/5}),
\end{equation*}
\begin{equation*}
  D(\eps,I)=\{\norm{Z}_{L^\infty(I)}<\eps\},
\end{equation*}
\begin{equation*}
  E(\eps,I)=\left\{\max_{\alpha=1,\ldots,n} \max_{i_1,\ldots,i_\alpha=1,\ldots,d} 
    \norm{A^{(\alpha)}_{i_1,\ldots, i_\alpha}}_{L^{\infty}(I)}< \eps\right\}.
\end{equation*}
Define 
\begin{equation*}
  \lambda^{(\alpha)}_{i_1,\ldots,i_\alpha}=A^{(\alpha)}_{i_1,\ldots,i_\alpha}(t_k).
\end{equation*}

On $G_k(\eps)$
\begin{align*}
  \|Z_\lambda\|_{L^\infty(I_k)} &\leq \|Z\|_{L^\infty(I_k)} +(n+1)d^n
  \max_{\alpha; i_1,\ldots,i_\alpha}
  \Lip(A^{(\alpha)}_{i_1,\ldots,i_\alpha})
  |t_{k}-t_{k-1}| \\
  & \qquad\ \qquad\ \qquad\ \qquad\ \qquad\ \qquad\ \qquad\ \qquad
  \max_{\alpha; i_1,\ldots,i_\alpha}\Norm{W_{i_1}\ldots
    W_{i_\alpha}}_{\frac14} \\ &\leq \eps^{8^{n+2}}
  +(n+1)d^n\eps^{-1}\eps^{8^{n+1}\cdot3/2}\eps^{-1/5}<\eps^{8^{n+1}\cdot5/4},
\end{align*}
for sufficiently small $\eps$, since
\begin{equation*}
  8^{n+2}>-1+8^{n+1}\cdot3/2-1/5>8^{n+1}\cdot5/4.  
\end{equation*}

On the other hand, for $\omega\in G_k$ there exists an $\alpha$ and
$i_1,\ldots,i_\alpha$ such that
\begin{equation*}
  \lambda^{(\alpha)}_{i_1,\ldots,i_\alpha}\geq\eps-
  \Lip(A^{(\alpha)}_{i_1,\ldots,i_\alpha})|t_{k}-t_{k-1}| \geq
  \eps-\eps^{-1}\eps^{8^{n+1}\cdot3/2}>\eps^{2+\frac{1}{n+1}}. 
\end{equation*}
Hence,
\begin{equation*}
  G_k(\eps)\subset
  D^*(\eps^{8^{n+1}\cdot5/4},I_k,\Lambda(\eps^{2+\frac{1}{n+1}},n))
  \subset F(\eps^{8^{n+1}\cdot5/4},\eps^{2+\frac{1}{n+1}}).    
\end{equation*}
\end{proof}

\begin{proof}[Proof of Theorem \ref{th:chance_of_being_small_estimate}]
  We begin by remarking that classical estimates on the supremum and
  H\"older continuity of a Wiener process combine to yield
  \begin{equation}\label{eq:weinerHolderSup}
    \PP(B^c(\eps^{-1/5}))<K_3e^{K_4\eps^{-q_2}},   
  \end{equation}
  for some positive $K_3(n),K_4(n),q_2(n)$.

  Theorem~\ref{th:chance_of_being_small_estimate} is then implied by
  the identity
  \begin{equation*}
    B^c(\eps^{-1/5})\cup F(\eps^{8^{n+1}\cdot5/4},\eps^{2+\frac{1}{n+1}})= 
    B^c(\eps^{-1/5})\cup
    \left( F(\eps^{8^{n+1}\cdot5/4},\eps^{2+\frac{1}{n+1}})\cap B(\eps^{-1/5}) \right),
  \end{equation*}
 the estimate from \eqref{eq:weinerHolderSup}, and the following lemma whose proof fills the remainder of this section.
\end{proof}

\begin{lemma}\label{lm:uniform_in_lambda} 
  For every $n$ there are positive numbers $q_3(n)$, $K_5(n)$,
  $K_6(n)$, and $\eps_2(n)$ such that for all $k=1,\ldots,n$,
  \begin{equation*}
    \PP(D^*(\eps^{8^{n+1}\cdot5/4},I_k,
    \Lambda(\eps^{2+\frac{1}{n+1}},n)) \cap B(\eps^{-1/5}))\leq K_5
    \exp\{-K_6\eps^{-{q_3}}\},  
  \end{equation*}
for $\eps<\eps_2$.
\end{lemma}

We shall derive this Lemma from the next one.
\begin{lemma} \label{lm:individual_lambda}For every $n$ there are
  positive numbers $q_4(n)$,  $K_7(n)$,  
$K_8(n)$, $\eps_3(n)$ with the
following property.

Let $\left\{\lambda^{(\alpha)}_{i_1,\ldots,
    i_\alpha},\alpha=0,\ldots,n, i_1,\ldots,
  i_\alpha=1,\ldots,d\right\}$ be a symmetric family of coefficients satisfying
\begin{equation*}
  \max \left\{|\lambda^{(n-\alpha)}_{i_1,\ldots, i_{n-\alpha}}|\eps^{-\left(2+\frac{1}{\alpha+1}\right)}:
    \alpha=0,\ldots,n, i_1,\ldots, i_\alpha=1,\ldots,d\right\}\geq1.
\end{equation*}
Define
\begin{equation*}
  D^*(\eps,I,\lambda)=\left\{\sup_{t\in I}|Z_\lambda(t)|<\eps\right\}.
\end{equation*}
Then
\begin{equation*}
  \PP(D^*(\eps^{8^{n+1}},I,\lambda)\cap B(\eps^{-1/5}))\leq
  K_7\exp\{-K_8\eps^{-{q_5}}\}, 
\end{equation*}
for $\eps<\eps_3$.
\end{lemma}

\bpf We shall prove this lemma by induction in $n$. If $n=0$, then the
statement of the lemma is obvious with the probability in the l.h.s. being
equal to $0$.

In the induction step we may always assume that
\begin{equation}
  \label{eq:normalizing_lambda_for_induction_step}
  \max \left\{|\lambda^{(n-\alpha)}_{i_1,\ldots, i_{n-\alpha}}|\eps^{-\left(2+\frac{1}{\alpha+1}\right)}:
    \alpha=0,\ldots,n, i_1,\ldots, i_\alpha=1,\ldots,d\right\}=1.
\end{equation}

Since the coefficients $\lambda$ are not random, we can use the It\^o formula
to write down the semimartigale representation of $Z_\lambda$, namely,
\begin{equation*}
  Z_\lambda(t)=V(t)+M(t),
\end{equation*}
where the finite variation part $V$ (which is, in fact, continuously
differentiable~a.s.) is given by

\begin{multline}
\label{eq:finite_variation_part}
  V(t)=\lambda^{(0)}+\sum_{i_1}\lambda_{i_1}^{(1)}W_{i_1}(t_1)+
         \sum_{i_1,i_2}\lambda_{i_1,i_2}^{(2)}W_{i_1}(t_1)W_{i_2}(t_1)\\+\ldots+
          \sum_{i_1,\ldots,i_n}\lambda^{(n)}_{i_1,\ldots,i_n}W_{i_1}(t_1)\ldots W_{i_n}(t_1)\\
   +\frac12\sum_{i_1,i_2} \lambda^{(2)}_{i_1i_2}\sum_{k_1\neq k_2}
   \int_{t_1}^{t}\frac{W_{i_1}(s)W_{i_2}(s)}{W_{i_{k_1}}(s)W_{i_{k_2}}(s)}\delta_{i_{k_1}i_{k_2}}ds
   \\+\ldots+\frac12\sum_{i_1,\ldots,i_n} \lambda^{(n)}_{i_1\ldots i_n}\sum_{k_1\neq k_2}
   \int_{t_1}^{t}\frac{W_{i_1}(s)\ldots W_{i_n}(s)}{W_{i_{k_1}}(s)W_{i_{k_2}}(s)}\delta_{i_{k_1}i_{k_2}}ds,
\end{multline}
and the martingale part $M$  is given by

\begin{multline*}
   M(t)=\sum_{i_1}\lambda^{(1)}_{i_1}\int_{t_1}^tdW_{i_1}(s)
   +\sum_{i_1,i_2} \lambda^{(2)}_{i_1i_2}\sum_{k}
   \int_{t_1}^{t}\frac{W_{i_1}(s)W_{i_2}(s)}{W_{i_{k}}(s)}dW_{i_{k}}(s)
   \\+\ldots+\sum_{i_1,\ldots,i_n} \lambda^{(n)}_{i_1\ldots i_n}\sum_{k}
   \int_{t_1}^{t}\frac{W_{i_1}(s)\ldots W_{i_n}(s)}{W_{i_{k}}(s)}dW_{i_{k}}(s).
\end{multline*}

For a function $f$ defined on a set $S$ denote
\begin{equation*}
  \osc_S f = \sup\{|f(s)-f(t)|:\ s,t\in S\}.
\end{equation*}

Since $\sup_I|Z|<\eps^{8^{n+1}}$ implies $\osc_I Z < 2\eps^{8^{n+1}}$,
the event of interest can be decomposed as
\begin{multline*}
  D^*(\eps^{8^{n+1}},I,\lambda)\cap B(\eps^{-1/5})\\
  \subset\left(\left\{\osc_I V<\eps^{8^{n+1}},
      \sup_{t\in I}|M(t)|<3\eps^{8^{n+1}}\right\}\cap B(\eps^{-1/5})\right)\\ 
  \cup\left(\{\osc_I V > \eps^{8^{n+1}}\}\cap B(\eps^{-1/5})\right).
\end{multline*}

For small $\eps$ the set $\{\osc_I V > \eps^{8^{n+1}}\}\cap B(\eps^{-1/5})$ 
in the decomposition above is empty. 
Indeed, \eqref{eq:normalizing_lambda_for_induction_step} implies that on this event 
each integral term with coefficient $\lambda_{i_1,\ldots,i_{n-\alpha}}^{(n-\alpha)}$ in~\eqref{eq:finite_variation_part} is bounded by
$\eps^{2+\frac{1}{\alpha+1}}\eps^{8^{n+1}\cdot3/2}\eps^{-1/5}<\eps^{8^{n+1}+\delta}$ for a positive $\delta$
and sufficiently small $\eps$, and there are only finitely many terms. 
Now,
\begin{multline*}
  \left\{\osc_I V<\eps^{8^{n+1}},\sup_{t\in I}|M(t)|<3\eps^{8^{n+1}}\right\}
  \cap B(\eps^{-1/5})
  \\
  \subset \left\{\sup_{t\in I}|M(t)|<3\eps^{8^{n+1}}\right\}\cap B(\eps^{-1/5})\\
  \subset  \left\{\sup_{t\in I}|M(t)|<3\eps^{8^{n+1}},
    \langle M\rangle_I>\eps^{8^{n+1}\cdot 15/8}\right\}\cap B(\eps^{-1/5})\\
  \cup \left\{\sup_{t\in I}|M(t)|<3\eps^{8^{n+1}},
    \langle M\rangle_I\leq\eps^{8^{n+1}\cdot 15/8}\right\}\cap B(\eps^{-1/5})\\
  \subset \left\{\sup_{t\in I}|M(t)|<3\eps^{8^{n+1}},
    \langle M\rangle_I>\eps^{8^{n+1}\cdot 15/8}\right\}
  \cup  \left(\left\{\langle M\rangle_I\leq\eps^{8^{n+1}\cdot 15/8}\right\}\cap B(\eps^{-1/5})\right).
\end{multline*}
Let us denote the sets in the r.h.s. by $D_1$ and $D_2$ respectively.
To estimate the probability of the set $D_1$ we need the following lemma
(see [Bass,p.209])
\begin{lemma} There exist $c_1,c_2>0$ such that if $M_t$ is a continuous
martingale, $T$ is a bounded stopping time, and $\eps>0$, then
\begin{equation*}
  \PP\left\{\sup_{t\leq T}|M_t|<\delta, \langle M\rangle_T>\eps\right\}\leq c_1e^{-c_2\eps/\delta^2}.
\end{equation*}
\end{lemma}

This result allows to conclude that
\begin{equation}
  \label{eq:F_is_small}
  \PP(D_1)\leq c_1\exp\left\{-\frac{c_2}{9}\eps^{-8^{n-1}}\right\}.
\end{equation}

To estimate $\PP(D_2)$, we notice that
the proof of Theorem \ref{th:quadratic_variation_of_polynomial} and the
continuous
differentiability of $V$ imply  that
\begin{equation*}
 \langle M\rangle_I= \sum_{r=1}^d\int_I\left(\sum_{\beta=1}^n
  \sum_{i_1,\ldots,i_\beta} 
  \sum_{p=1}^\beta\delta_{ri_p} 
  \frac{\lambda^{(\beta)}_{i_1,\ldots,i_\beta}W_{i_1}
    \ldots W_{i_\beta}}{W_{i_p}}\right)^2 ds.
\end{equation*}
Therefore,
\begin{equation*}
  \PP(D_2)\leq \min_{r=1,\ldots,d}
  \PP\left(D_2(r)\cap  B(\eps^{-1/5})\right),
\end{equation*}
where
\begin{equation*}
  D_2(r) = \left\{\int_I\left(\sum_{\beta=1}^n
      \sum_{i_1,\ldots,i_\beta} \sum_{p=1}^\beta\delta_{ri_p} \frac{\lambda^{(\beta)}_{i_1,\ldots,i_\beta}W_{i_1}
        \ldots W_{i_\beta}}{W_{i_p}}\right)^2 ds<\eps^{8^{n+1}\cdot15/8}\right\}.
\end{equation*}

There exist $\beta$ and $i_1\ldots i_{\beta}$ such that
$|\lambda^{(\beta)}_{i_1\ldots
  i_{\beta}}|=\eps^{2+\frac{1}{n-\beta+1}}$.  If $\beta\neq 0$, then
choose $r$ so that the definition of $D_2(r)$ contains that
$\lambda^{(\beta)}_{i_1\ldots i_{\beta}}$ and define
\begin{equation*}
  Z_{\lambda,r}=\sum_{\alpha=1}^n  \sum_{i_1,\ldots,i_\alpha} 
  \sum_{p=1}^\alpha\delta_{ri_p}   \frac{\lambda^{(\alpha)}_{i_1,\ldots,i_\alpha}W_{i_1}
    \ldots W_{i_\alpha}}{W_{i_p}}.
\end{equation*}

We want to prove that
\begin{equation}
  \label{eq:induction_inclusion}
  D_2(r)\cap  B(\eps^{-1/5})\subset \left\{\sup_{t\in I}|Z_{\lambda,r}(t)|
    <\eps^{8^n}\right\}\cap B(\eps^{-1/5}).
\end{equation}

On the set $B(\eps^{-1/5})$ the H\"older constant of $Z_{\lambda,r}$
is bounded by $ nd^n\eps^{-1/5}\eps^{2+1/(n+1)}$. So,
if the condition $\sup_{t\in I}|Z_{\lambda,r}(t)|<\eps^{8^n}$ is not
fullfilled, we have
\begin{equation*}
  \inf_{t\in I}|Z_{\lambda,r}(t)|\geq \eps^{8^n}-nd^n\eps^{-1/5}\eps^{2+1/(n+1)}
  \left(\eps^{8^{n+1}\cdot3/2}\right)^{1/4}\geq c\eps^{8^n},
\end{equation*}
for some constant $c$, since
$8^n<-1/5+2+1/(n+1)+8^{n+1}\cdot3/8$.
Thus, on $D_2(r)$ we have
\begin{equation*}
  \eps^{8^{n+1}\cdot15/8}>\int_IZ_{\lambda,r}^2(s)ds\geq\frac12\eps^{8^{n+1}\cdot3/2}
  (c\eps^{8^n})^2=\frac{c^2}{2}\eps^{8^{n+1}\cdot{14/8}},
\end{equation*}
which is impossible for small $\eps$. Therefore, our assumption was false and
\eqref{eq:induction_inclusion} is proved. Now \eqref{eq:induction_inclusion}
and the induction assumption imply
\begin{multline}
  \label{eq:induction_inclusion_estimate}
  \PP\left(D_2(r)\cap  B(\eps^{-1/5})\right)\leq \PP\left(\left\{\sup_{t\in I}|Z_{\lambda,r}(t)|
      <\eps^{8^n}\right\}\cap B(\eps^{-1/5})\right)\\\leq K_7(n-1)\exp\{-K_8(n-1)\eps^{-{q_5(n-1)}}\}.
\end{multline}

Consider now the case where  $|\lambda^{(\beta)}_{i_1\ldots i_\beta}|<\eps^{2+\frac{1}{n-\beta+1}}$
if $\beta \neq 0,$ and $|\lambda^{(0)}|=\eps^{2+\frac{1}{n+1}}$. 
Denote 
\begin{equation*}
  W^*=\sup\{|W_{i_1}\ldots W_{i_\alpha}|:
  \alpha=1,\ldots,n-1,i_1,\ldots,i_\alpha=1,\ldots,d \}. 
\end{equation*}
We have
\begin{equation*}
  \PP(D(\eps^{8^{n+1}},I,\lambda))\leq\PP\{\eps^{2+\frac{1}{n+1}}-nd^n\eps^{2+\frac{1}{n}}
  W^*<\eps^{8^{n+1}}\}
  = \PP\left\{W^*>\frac{\eps^{2+\frac{1}{n+1}}-\eps^{8^{n+1}}}{nd^n\eps^{2+\frac{1}{n}}}\right\}.
\end{equation*}
Since $8^{n+1}>2+\frac{1}{n+1}$ and $2+\frac{1}{n}>2+\frac{1}{n+1}$,
\begin{equation*}
  \PP(D(\eps^{8^{n+1}},I,\lambda))\leq K_9\exp\{K_{10}(n)\eps^{-q_{6}}\},
\end{equation*}
for some positive constants $K_9(n),K_{10}(n),q_{6}(n)$.  This
completes the proof of the lemma.\epf

\bpf[Proof of Lemma~\ref{lm:uniform_in_lambda}] It suffices to show
\begin{equation*}
  \PP(D^*(\eps^{8^{n+1}\cdot5/4},I,
  \overline\Lambda(\eps^{2+\frac{1}{n+1}},n)) \cap B(\eps^{-1/5}))\leq
  K_{11} \exp\{-K_{12}\eps^{-{q_7}}\},  
\end{equation*}
for some positive constants $K_{11}(n),K_{12}(n),q_{7}(n)$,
where
$\overline\Lambda(\eps,n)$ is the set of vectors  
$\left\{\lambda^{(\alpha)}_{i_1,\ldots, i_\alpha},\alpha=0,\ldots,n,
  i_1,\ldots, i_\alpha=1,\ldots,d\right\}$ 
such that
\begin{equation*}
  \max \left\{|\lambda^{(\alpha)}_{i_1,\ldots,
      i_{\alpha}}|:\alpha=0,\ldots,n, i_1,\ldots,
    i_\alpha=1,\ldots,d\right\}=\eps. 
\end{equation*}

For sufficiently small $\delta>0$ there is a set of points
$\{\lambda(\delta,j),
j=1,\ldots,[\delta^{-(n+1)d^n}]\}\subset\overline\Lambda(\eps^{2+\frac{1}{n+1}},n)$
such that for every
$\lambda\in\overline\Lambda(\eps^{2+\frac{1}{n+1}},n)$ there is $j$
such that $|\lambda(\delta,j)^{(\alpha)}_{i_1,\ldots,i_\alpha}-
\lambda^{(\alpha)}_{i_1,\ldots,i_\alpha}| <
\eps^{2+\frac{1}{n+1}}\delta$ for all $\alpha$ and
$i_1,\ldots,i_\alpha$. This implies
\begin{equation*}
  |Z_{\lambda{(\delta,j)}}-Z_\lambda|\leq
  (n+1)d^n\eps^{2+\frac{1}{n+1}}\delta\eps^{-1/5}. 
\end{equation*}
Choose $\delta=\eps^{8^{n+1}\cdot3/2}$. If $\sup_{t\in
  I}|Z_\lambda|<\eps^{8^{n+1}\cdot5/4}$, then
\begin{equation*}
  |Z_{\lambda(\eps^{8^{n+1}\cdot3/2},j)}|\leq \eps^{8^{n+1}\cdot5/4}+(n+1)d^n
  \eps^{2+\frac{1}{n+1}}\eps^{8^{n+1}\cdot3/2}\eps^{-1/5}<\eps^{8^{n+1}}.
\end{equation*}

Therefore, Lemma~\ref{lm:individual_lambda} implies
\begin{multline*}
  \PP\left(D^*(\eps^{8^{n+1}\cdot5/4},I,\overline\Lambda(\eps^{2+\frac{1}{n+1}},n))\cap
    B(\eps^{-1/5})\right)  \\
  \leq[\delta^{-(n+1)d^n}+1]
  \sup_{\lambda\in\overline\Lambda(\eps^{2+\frac{1}{n+1}},n) }
  \PP\left(D^*(\eps^{8^{n+1}},I,\lambda) \cap B(\eps^{-1/5})\right)
  \\\leq K_7\exp\{-K_8\eps^{-{q_5}}\}.
\end{multline*}
This completes the proof of Lemma~\ref{lm:uniform_in_lambda}.
\epf

\section {Polynomial vector fields. Derivatives and Lie brackets}
\label{sc:poynomial}

We start with a characterization of multilinear continuous operators,
which is an obvious generalization of the linear case:
\begin{lemma} \label{l:continuous_Poly_in_terms_of_norms} Let $\XX$ and $\YY$ be two
  Banach spaces.  Let $Q:\XX^{m} \to \YY$ be
  an $m$-linear operator which is continuous at zero. Then 
  \begin{equation*}
    |Q(x_1,\ldots,x_m)|_{\YY} \leq c |x_1|_{\XX}\cdots|x_m|_{\XX},
  \end{equation*}
 where 
  \begin{equation*}
    c=\sup_{|x_1|_{\XX},\ldots,|x_m|_{\XX}\le 1} 
|Q(x_1,\ldots,x_m)|_{\YY}.
  \end{equation*}
\end{lemma}

We define the local Lipschitz constant for a map $Q:\XX \to \YY$
as
\begin{equation}\label{eq:lipDef}
    \Lip(Q)(x)=\lim_{\eps \to 0} \sup_{\substack{z\in \XX\\  |z|_\XX
    \leq \eps}} \frac{|Q(x)-Q(z)|_\YY}{|x-z|_\XX}\;.
\end{equation}

\begin{lemma}\label{l:locallipPoly} Let $\XX$ and $\YY$ be two
  Banach spaces. Suppose $Q:\XX \to \YY$ is a continuous polynomial vector field of
  order $m$.
Then there is a constant $c$ such that
\begin{equation*}
   \Lip(Q)(x) \leq c(1+|x|_\XX)^m. 
\end{equation*}
\end{lemma}
\begin{proof}
This is an easy consequence of
Lemma~\ref{l:continuous_Poly_in_terms_of_norms}, since the
latter
implies a straightforward bound on the local Lipschitz constant for 
$Q$ in each of the $m$ variables.
\end{proof}

The  Fr\'echet derivative of order $i$ of a function $Q:\VV\to\VV'$ at a
point $y$ will be denoted by $(\FD^iQ)(y):\VV^i\to\VV'$. It is an
$i$-linear operator and its value at a tangent
vector~$(\phi_1,\ldots,\phi_i)\in\VV^i$ is denoted by $(\FD^iQ)(y)(\phi_1,\ldots,\phi_i)$.

\begin{lemma}\label{lm:higher_derivatives_of_multilinear}
Let $Q$ be a $j$-linear symmetric function. Then,
\begin{equation*}
  (\FD^i Q)(y)(\psi_1,\ldots,\psi_n)= \left\{\begin{array}{ll}\frac{j!}{(j-i)!}
    Q(y^{\otimes(j-i)},\psi_1,\ldots,\psi_i),&i\leq j,\\ \\
    0,& i>j.\end{array} \right.
\end{equation*}
\end{lemma}
\bpf If $i=1$, the Lemma immediately follows from the chain rule.  The
general case follows from an iterative application of the statement
for $i=1$. 
\epf

\begin{lemma}\label{lm:bound_higher-order_derivavtives} If
  $Q:\VV\to\VV'$ is a polynomial vector field of order $m$ such that
  condition~\eqref{eq:Poly1} holds true, then for every
  $i=2,\ldots,m$ there is a constant $K_i>0$ such that
  \begin{equation*}
    |D^{(i)}Q(y)(\psi_1,\ldots, \psi_i)|\leq K_i(1+\|y\|^{m-i})\|\psi_1\|\ldots\|\psi_i\|.
  \end{equation*}
\end{lemma}
\begin{proof}
   This is a straightforward consequence of Lemma~\ref{lm:higher_derivatives_of_multilinear}.
\end{proof}

\begin{lemma}\label{lm:brackets_with_constant_fields}
  Suppose $f_1,\ldots,f_i\in\VV$ are constant vector fields and
  $Q(x):\VV\to\VV'$ is a Fr\'echet differentiable vector field. Then
\begin{equation*}
  (D^iQ)(x)(f_1,\ldots,f_i)=[Q,f_1,f_2,\ldots,f_i](x).
\end{equation*}
\end{lemma}

\begin{proof}
  The lemma is proved by induction:
\begin{align*}
  (D^iQ)(x)(f_1,\ldots,f_i)&=(D(D^{i-1}Q)(\cdot)(f_1,\ldots,f_{i-1}))(x)(f_i)\\
  &=[D^{i-1}Q(\cdot)(f_1,\ldots,f_{i-1}),f_i](x).  
\end{align*}
\end{proof}

Recall that 
\begin{equation*}
  X(t)=u(0)+\int_0^tF(u(s))ds+\int_0^tf(s)ds.
\end{equation*}

\begin{lemma}\label{l:throughVF}\label{lm:expanding_polylynear} 
  Let $Q:\VV \to \VV'$ be a polynomial vector field. Then
  \begin{align*}
    Q\big(u(t)\big) =Q(X(t))+
    \sum_{i=1}^{n}\sum_{k_1,\ldots,k_i}[Q,g_{k_1},\ldots,g_{k_i}](X(t))\;W_{k_1}\ldots
    W_{k_i}.
  \end{align*}
\end{lemma}
\begin{proof}
    Since $Q$ is polynomial, we have
  \begin{equation*}
    Q(y)=\sum_{j=0}^n Q_j(y^{\otimes j})    
  \end{equation*}
  for some $n\in\N$ where $Q_j$ is a continuous, symmetric, multilinear
  vector field for each $j$. 
Now
  \begin{align*}
    Q(u(t))&=\sum_{j=0}^n Q_j((X(s)+\sum_kg_kW_k)^{\otimes j})\\
    &=\sum_{j=0}^n\sum_{i=0}^{j}\sum_{k_1,\ldots,k_i}
    \frac{j!}{(j-i)!}Q_j(X(s)^{\otimes(j-i)},g_{k_1},\ldots,g_{k_i})W_{k_1}\ldots
    W_{k_i}
    \\&= \sum_{i=0}^{n}\sum_{k_1,\ldots,k_i}\sum_{j=i}^n
    \frac{j!}{(j-i)!}Q_j(X(s)^{\otimes(j-i)},g_{k_1},\ldots,g_{k_i})W_{k_1}\ldots
    W_{k_i}.
  \end{align*}
Using Lemma \ref{lm:higher_derivatives_of_multilinear} 
we  have 
\begin{align*}
  \sum_{j=i}^n
    \frac{j!}{(j-i)!}Q_j(X(s)^{\otimes(j-i)},g_{k_1},\ldots,g_{k_i})&=  \sum_{j=0}^n
    [Q_j,g_{k_1},\ldots,g_{k_i}](X(s)) \\
    &=[Q,g_{k_1},\ldots,g_{k_i}](X(s)), 
\end{align*}
which completes the proof.
\end{proof}

\section{Bounds on norms and Lipschitz constants}
\label{sc:bounds_on_lipschitz_norms}
We define
\begin{equation*} 
\LipT \|K_{s,t}\|_{\VV \to \VV'}=\sup_{\|\phi\|\le1}\LipT \|K_{s,t}\phi\|_{\VV'}
\end{equation*}
and
\begin{equation*} 
\Norm{K_{s,t}}_\VV=\sup_{\|\phi\|\le 1}\Norm{K_{s,t}\phi}_{\VV}.
\end{equation*}
\begin{lemma}\label{l:lipBoundsXK}
  Under Assumptions \ref{a:solReg} and \ref{a:K2a}, for any $p \geq 1$,
  there is a universal constant $C_p$ such that the following
  bounds hold:
  \begin{align*}
    \E \big(\;\LipT \|K_{s,t}\|_{\VV \to \VV'}\big)^p &\leq C_p
    \sqrt{K^*_{2p}(1 + u^*_{2p(m-1)})} \\
    \E \Norm{K_{s,t}}_\VV &\leq C_p \sqrt{K^*_{2p}(1 + u^*_{2p(m-1)})} +
    K_{2p}^*.
\end{align*}
\end{lemma}
\begin{proof}
  From the equation for $K_{s,t}$ and the bound
  on $F$ (and hence $DF*$) from \eqref{eq:Poly1}, we see that
  for $\phi \in \VV$ with $ \|\phi\| \le1$
  \begin{align*}
    \LipT \|K_{s,t}\phi\|_{\VV'} \leq&
    \supT\|DF^*(u(r))K_{r,t}\phi\|_{\VV'}\\
    \leq & C( 1 +
    \supT\|u(s)\|^{m-1})\;\supT\|K_{s,t}\phi\|.
  \end{align*}
  Next, take the supremum over $\phi$, then the $p$th power, and lastly the
  expected value.  The first inequality of the Lemma follows from the bounds in 
  Assumptions~\ref{a:solReg} and~\ref{a:K2a} after applying the Cauchy-Schwartz
  inequality to the right hand side. The second inequality 
  follows from the first one
  and the assumptions. 
\end{proof}

\begin{lemma}\label{l:lipPoly}
  Let $Q:\HH \to \VV'$ be a continuous polynomial vector field of
  order $m$ and let $f_i:[0,T] \to \HH$ for $i \in
  \{1,\cdots,m\}$. Then there exists  a universal  constant $c(m)$
  such that
  \begin{align*}
    \LipT_{\VV'}\Big(Q(f_1(t),\cdots,f_m(t))\Big) \leq c \prod_{i=1}^{m}
    \big(1+\supT_{\HH}(f_i)\big) \sum_{i=1}^m \LipT_{\HH}(f_i).
  \end{align*}
\end{lemma}

\begin{proof}
The proof is analogous to that of Lemma~\ref{l:locallipPoly}
\end{proof}

\begin{lemma}\label{l:boundLipSupIP} If $f,g:[0,T] \to \VV$  then
  \begin{align*}
    \LipT|\ip{f}{g}| & \leq \LipT|f|_{\VV'}\supT\|g\|+ \LipT|g|\supT|f|\\
    \Norm{\ip{f}{g}}& \leq {\Norm{f}^2_\HH + \Norm{g}^2_\VV}.
  \end{align*}
  \end{lemma}
  \begin{proof}
    The first bound follows from
    \begin{align*}
      |\ip{f(t)}{g(t)}-\ip{f(s)}{g(s)}| &\leq |\ip{f(t)-f(s)}{g(t)}|
      +|\ip{f(s)}{g(t)-g(s)}|\\
    &\leq  |f(t)-f(s)|_{\VV'} \|g(t)\|  +     |f(t)| |g(t)-g(s)|.
    \end{align*}

    We turn to the second bound. Since 
    \begin{align*}
      \supT|\ip{f}{g}| &\leq \supT|f|\; \supT|g| \leq \frac{\supT
        |f|^2 + \supT\|g\|^2}{2}\\
      &\leq  {\Norm{f}_\HH^2 + \Norm{g}_\VV^2},
    \end{align*}
    the first inequality of the lemma implies
    \begin{align*}
       \LipT|\ip{f}{g}| \leq  {\Norm{f}_\HH^2 + \Norm{g}_\VV^2},
    \end{align*}
    and we are done.
  \end{proof}

\bibliographystyle{alpha}

\begin{thebibliography}{BAL91}

\bibitem[AKSS]{AgrachevKuksinSarychevShirikyan06FDP}
A.~Agrachev, S.~Kuksin, A.~Sarychev, and A~Shirikyan.
\newblock On finite-dimensional projections of distributions for solutions of
  randomly forced pde's.
\newblock Preprint 2006.

\bibitem[BAL91]{b:BenArousLeandre91}
G.~Ben~Arous and R.~L{\'e}andre.
\newblock D\'ecroissance exponentielle du noyau de la chaleur sur la diagonale.
  {I}.
\newblock {\em Probab. Theory Related Fields}, 90(2):175--202, 1991.

\bibitem[Bas98]{bass99DEO}
Richard~F. Bass.
\newblock {\em Diffusions and elliptic operators}.
\newblock Probability and its Applications (New York). Springer-Verlag, New
  York, 1998.

\bibitem[Bel87]{b:Bell87}
Denis~R. Bell.
\newblock {\em The {M}alliavin calculus}.
\newblock Longman Scientific \& Technical, Harlow, 1987.

\bibitem[BT05]{BaudoinTeichmann05HID}
Fabrice Baudoin and Josef Teichmann.
\newblock Hypoellipticity in infinite dimensions and an application in interest
  rate theory.
\newblock {\em Ann. Appl. Probab.}, 15(3):1765--1777, 2005.

\bibitem[Cer99]{Cerrai99SPT}
Sandra Cerrai.
\newblock Smoothing properties of transition semigroups relative to {SDE}s with
  values in {B}anach spaces.
\newblock {\em Probab. Theory Related Fields}, 113(1):85--114, 1999.

\bibitem[CF88]{b:CoFo88}
Peter Constantin and Ciprian Foia{\c{s}}.
\newblock {\em {N}avier-{S}tokes Equations}.
\newblock University of Chicago Press, Chicago, 1988.

\bibitem[DL88]{b:DautrayLions88}
Robert Dautray and Jacques-Louis Lions.
\newblock {\em Analyse math\'ematique et calcul num\'erique pour les sciences
  et les techniques.}
\newblock Masson, Paris, 1988.

\bibitem[DPZ96]{b:DaZa96}
Giuseppe Da~Prato and Jerzy Zabczyk.
\newblock {\em Ergodicity for Infinite Dimensional Systems}.
\newblock Cambridge, 1996.

\bibitem[EH01]{EckmannHairer00UIM}
J.-P. Eckmann and M.~Hairer.
\newblock Uniqueness of the invariant measure for a stochastic {PDE} driven by
  degenerate noise.
\newblock {\em Comm. Math. Phys.}, 219(3):523--565, 2001.

\bibitem[EM01]{short:EMattingly01}
Weinan E and Jonathan~C. Mattingly.
\newblock Ergodicity for the {N}avier-{S}tokes equation with degenerate random
  forcing: finite-dimensional approximation.
\newblock {\em Comm. Pure Appl. Math.}, 54(11):1386--1402, 2001.

\bibitem[Fla94]{b:Fl94}
Franco Flandoli.
\newblock Dissipativity and invariant measures for stochastic {N}avier-{S}tokes
  equations.
\newblock {\em NoDEA}, 1:403--426, 1994.

\bibitem[HM06]{b:HairerMattingly04}
Martin Hairer and Jonathan~C. Mattingly.
\newblock Ergodicity of the degenerate stochstic {2D} {N}avier--{S}tokes
  equation.
\newblock {\em Annals of Mathematics}, 164(3), 2006.

\bibitem[Kli87]{b:Kliemann87}
Wolfgang Kliemann.
\newblock Recurrence and invariant measures for degenerate diffusions.
\newblock {\em The Annals of Probability}, 2:690--707, 1987.

\bibitem[KS84]{b:KusuokaStroock84}
Shigeo Kusuoka and Daniel Stroock.
\newblock Applications of the {M}alliavin calculus. {I}.
\newblock In {\em Stochastic analysis (Katata/Kyoto, 1982)}, pages 271--306.
  North-Holland, Amsterdam, 1984.

\bibitem[Mal78]{b:Malliavin78}
Paul Malliavin.
\newblock Stochastic calculus of variation and hypoelliptic operators.
\newblock In {\em Proceedings of the International Symposium on Stochastic
  Differential Equations (Res. Inst. Math. Sci., Kyoto Univ., Kyoto, 1976)},
  pages 195--263, New York, 1978. Wiley.

\bibitem[MB02]{b:MajdaBertozzi02}
Andrew~J. Majda and Andrea~L. Bertozzi.
\newblock {\em Vorticity and incompressible flow}, volume~27 of {\em Cambridge
  Texts in Applied Mathematics}.
\newblock Cambridge University Press, Cambridge, 2002.

\bibitem[MP06]{MattinglyPardoux05}
Jonathan~C. Mattingly and {\'E}tienne Pardoux.
\newblock Malliavin calculus for the stochastic 2{D} {N}avier-{S}tokes
  equation.
\newblock {\em Comm. Pure Appl. Math.}, 59(12):1742--1790, 2006.

\bibitem[Nor86]{b:Norris86}
James Norris.
\newblock Simplified {M}alliavin calculus.
\newblock In {\em S\'eminaire de Probabilit\'es, XX, 1984/85}, pages 101--130.
  Springer, Berlin, 1986.

\bibitem[Nua95]{b:Nualart95}
David Nualart.
\newblock {\em The {M}alliavin calculus and related topics}.
\newblock Probability and its Applications. Springer-Verlag, New York, 1995.

\bibitem[Oco88]{b:Ocone88}
Daniel Ocone.
\newblock Stochastic calculus of variations for stochastic partial differential
  equations.
\newblock {\em J. Funct. Anal.}, 79(2):288--331, 1988.

\bibitem[Rom04]{b:Romito04}
Marco Romito.
\newblock Ergodicity of the finite dimensional approximation of the 3{D}
  {N}avier-{S}tokes equations forced by a degenerate noise.
\newblock {\em J. Statist. Phys.}, 114(1-2):155--177, 2004.

\bibitem[SY02]{SellYou02DEE}
George~R. Sell and Yuncheng You.
\newblock {\em Dynamics of evolutionary equations}, volume 143 of {\em Applied
  Mathematical Sciences}.
\newblock Springer-Verlag, New York, 2002.

\bibitem[Wu06]{Wu2006Thesis}
Ming-Yih Wu.
\newblock {\em Stochastic Boussinesq Equations and the infinite dimensional
  Malliavin Calculus}.
\newblock PhD thesis, Princeton, 2006.

\end{thebibliography}

\def\Rom#1{\uppercase\expandafter{\romannumeral #1}}\def\cprime{$'$}

\end{document}